\documentclass[a4paper,fleqn]{cas-sc}

\usepackage[numbers,sort&compress]{natbib}
\usepackage{subfigure}
\graphicspath{{pic/}}

\def\tsc#1{\csdef{#1}{\textsc{\lowercase{#1}}\xspace}}
\tsc{WGM}
\tsc{QE}
\tsc{EP}
\tsc{PMS}
\tsc{BEC}
\tsc{DE}

\begin{document}
\let\WriteBookmarks\relax
\def\floatpagepagefraction{1}
\def\textpagefraction{.001}
\shorttitle{Physics-informed Multiresolution Wavelet Neural Network Method for Solving Partial Differential Equations}
\shortauthors{F. Han et~al.}

\title [mode = title]{Physics-informed Multiresolution Wavelet Neural Network Method for Solving Partial Differential Equations}                      
\tnotemark[1]

\tnotetext[1]{This document is the results of the research
   project funded by the National Key Research and Development Program of China 2020YFA0709800 and National Natural Science Foundation of China (Grant No. 92470115).}

\author[1]{Feng Han}[style=chinese,
                        orcid=0000-0002-6729-9778]

\credit{Conceptualization of this study, Methodology, Writing - Original draft preparation}

\author[1,2]{Jianguo Wang}[style=chinese,
                        orcid=0000-0002-9457-9425]
\cormark[1]
\ead{wanguiuc@mail.xjtu.edu.cn}
\credit{Conceptualization of this study, Methodology, and Checking draft preparation}

\author[1]{Guoliang Peng}[style=chinese,
                        orcid=0009-0006-5131-8694]
\credit{Data curation, Implemented the algorithm}

\author[1]{Xueting Shi}[style=chinese,
                        orcid=0009-0001-8476-2058]
\credit{Data curation, Implemented the code of test cases}

\affiliation[1]{organization={National Key Laboratory of Intense Pulsed Radiation Simulation and Effect; Northwest Institute of Nuclear Technology},
                city={Xi'an},
                postcode={710024}, 
                state={Shaanxi},
                country={China}}

\affiliation[2]{organization={School of Information and Communications Engineering, Xi’an Jiaotong University},
                city={Xi'an},
                postcode={710049}, 
                state={Shaanxi},
                country={China}}

\cortext[cor1]{Corresponding author: Jianguo Wang}

\begin{abstract}
In this paper, a physics-informed multiresolution wavelet neural network (PIMWNN) method is proposed for solving partial differential equations (PDEs). This method uses the multiresolution wavelet neural network (MWNN) to approximate unknown functions, then substituting the MWNN into PDEs and training the MWNN by least-squares algorithm. We apply the proposed method to various problems, including stationary/nonstationary advection, diffusion and advection-diffusion problems, and linear/nonlinear time-dependent problems. Numerical experiments show that the PIMWNN method can achieve higher accuracy and faster speed than Physics Informed Neural Networks (PINNs). Moreover, the PIMWNN method, being mesh-free, can handle different boundary conditions easily and solve the time-dependent problems efficiently. The proposed method is expected to solve the spectral bias problem in network training. These characteristics show the great potential of the PIMWNN method used in the field of numerical solving methods for PDEs.
\end{abstract}


\maketitle

\section{Introduction}

Artificial neural networks (ANNs), also known as neural networks, have found extensive applications in traditional machine learning tasks, including image classification, image generation, and complex game playing \cite{EWN2020,Karniadakis2021PhysicsinformedML}. Over the past decade, deep learning, particularly through deep neural networks, has been widely adopted across various disciplines. Recently, the use of deep learning to solve partial differential equations (PDEs) that govern physical phenomena has emerged as a vibrant area of research in scientific computation \cite{raissi2019,raissi2020}. In \cite{raissi2019}, Raissi et al. introduced the physics-informed neural networks (PINNs) algorithm. This approach embeds the strong form of the governing equations into the neural network and constructs the loss function by incorporating a residual term derived from the equations. This residual term acts as a penalty to restrict the solution space, effectively transforming the PDE problem into an optimization task focused on minimizing the loss function. PINNs and their variants have been successfully applied to both forward and inverse problems involving PDEs, including nonlinear conservation laws, the Navier-Stokes equations, the Korteweg-de Vries (KdV) equations, compressible Euler equations \cite{raissi2019,jagtap2020}, stochastic PDEs \cite{yangliu2020,dongkunzhang2019,dongkunzhang2020,YiboYang2019}, fractional PDEs \cite{Pang2019}, and so on.

Alternative approaches to solve PDEs using deep learning leverage the weak form of the equations. For example, the deep Ritz method employs the variational formulation of symmetric PDEs and defines the loss function through the integral of the boundary residual. This method trains deep neural networks to minimize the numerical quadrature of the loss function \cite{EWN2018}. Extending this framework, the Deep Nitsche method (DNM) was developed to handle essential boundary conditions within the deep Ritz method \cite{Ming2021}. Additionally, various neural network architectures based on the Petrov-Galerkin scheme have been proposed, offering diverse choices for trial functions and test spaces \cite{Khodayi2020,Yang2021hpVariationalPN,Berrone2021,Kharazmi2019VariationalPN}. Similarly, DeepONets \cite{swang2021} and neural operators \cite{Kovachki2023} provide alternative ansatz spaces, broadening the scope of physics-informed machine learning models.

Despite the success of these methods, several challenges persist. The optimization of neural network objective functions during the solving process is highly dependent on the characteristics of the optimization problems and algorithms, leading to issues such as limited computational accuracy, prolonged training times, and difficulties in capturing high-frequency components of solutions. Furthermore, there is a lack of a robust theoretical foundation for selecting neural networks that approximate PDE solutions, often relying on empirical experience \cite{EWN2020,raissi2019,Khodayi2020,Yang2021hpVariationalPN,Berrone2021,Kharazmi2019VariationalPN,Uddin2023WaveletsBP,math10121976}. To address some of these challenges, \cite{Wangfei2023} proposed a PDE-solving method based on random weighted neural networks combined with the least squares method, which improves accuracy and reduces computational costs. Then, this method is combined with the discontinuous finite element method and the finite difference method to solve problems such as the linear elasticity and Navier-Stokes equations, the wave equation, and the Stokes equation \cite{YongShang2024,SUN2024128,SUN2024115830,DangHaoning2024,LI2025113847}. Additionally, the physics-informed extreme learning machine (PIELM), a faster variant of PINNs, was introduced for solving stationary and time-dependent linear PDEs \cite{Dwivedi2019PhysicsIE}.

Wavelet neural networks (WNNs) represent another innovative approach, combining wavelet theory with neural network theory \cite{Zhang1992WaveletN,Arts2022TheFC}. One type of WNN integrates wavelet bases with feed-forward neural networks, using wavelet functions as activation functions in the hidden layer. Due to the localized properties of wavelet activation functions, WNNs offer advantages in generalization capability and learning speed compared to traditional networks \cite{Zainuddin2011}. WNNs have been successfully applied in various real-world applications, such as pattern classification, forecasting, and function approximation \cite{Rafiei2017}. In \cite{Uddin2023WaveletsBP}, wavelet neural networks were used to replace fully connected feed-forward neural networks in the PINN method, resulting in a wavelet-based PINN for solving nonlinear differential equations. Another advancement is the deep wavelet neural network, which incorporates a wavelet mapping layer into the fully connected feed-forward neural network of the PINN method \cite{math10121976}. A further development is the multiresolution wavelet neural network (MWNN), which combines multiresolution wavelet analysis with neural networks \cite{Bakshi1993}. The MWNN is a feed-forward neural network with a single hidden layer, where the activation functions are a family of orthonormal wavelets. According to multiresolution analysis theory, the MWNN can approximate any function in the square integrable function space, with its approximation capability improving as the number of hidden layer nodes increases.

In this paper, we propose the physics-informed multiresolution wavelet neural network (PIMWNN) method for solving PDEs. The method involves three key steps: first, approximating numerical solutions using the multiresolution wavelet neural network; second, substituting the MWNN into the PDEs and initial/boundary conditions to derive a system of linear equations, where the unknowns are the weight coefficients of the MWNN; and third, training the MWNN by solving a least-squares problem associated with the linear system. This approach aims to combine the strengths of physics-informed learning with the multiresolution capabilities of wavelet neural networks, offering a promising framework for efficient and accurate PDE solutions.

The main contributions of this work can be listed as follows,

\begin{itemize}
\item We propose a physics-informed multiresolution wavelet neural network method, which makes the network structure interpretable. With the increase of nodes in the hidden layer of the network, the approximation ability of the network is increasing.
\item Compared with the training method based on an optimization of the loss function, the PIMWNN training is completed by solving linear equations, which greatly improves the solving speed of problems.
\item The PIMWNN method can capture the high-frequency part of the solutions and solve the spectral bias problem in the network training with the increase of the number of wavelet scale. 
\end{itemize}

The rest of the paper is organized as follows: In Section \ref{sec:second}, we review some preliminaries of the problem. Section \ref{method} gives the methodology of the proposed method including a brief discussion on structure of the neural network, PIMWNN methods for linear and nonlinear problems and so on. Numerical experiments are performed to show the effectiveness of PIMWNN methods in Section \ref{num_exp}. Finally, we end the paper with a summary and discussion in Section \ref{discussion}.

\section{Preliminaries}\label{sec:second}
\subsection{Brief review of the multiresolution analysis}\label{sec:model}
A multiresolution approximation is a sequence of nested space $V_j\subset V_{j+1}\subset\cdots\subset L^2(\mathbb{R}),j\in \mathbb{Z}$, such that:
\begin{itemize}
\item $\bigcap{V_j} = \{0\}, \overline{\bigcup V_j} = L^2$,
\item $f(t)\in V_j \Leftrightarrow f(2t)\in V_{j+1}$.
\end{itemize}
There exists a scaling function $\phi \in V_0$ such that $\phi_{j,k}(t) = 2^{j/2}\phi(2^jt-k),k\in \mathbb{Z},$ constitute a Riesz basis of $V_j.$

The scaling function $\phi\in V_0\subset V_1$ should satisfy a two-scale equation
\begin{equation}
    \phi(t) = \sum_{k\in\mathbb{Z}}h_k\phi(2t-k),
\end{equation}
where $h_k=\int \phi(t)\phi(2t-k)dt$.

Assuming that $\{\phi_{j,k}\}_{k\in\mathbb{Z}}$ is an orthonormal basis, the wavelet $\psi(t)$ is constructed as
\begin{equation}
    \psi(t) = \sum_{k\in\mathbb{Z}}g_k\phi(2t-k),
\end{equation}
with $g_k=(-1)^k h_{1-k}$, then $\{\psi_{j,k}\}_{k\in\mathbb{Z}}$ is an orthonormal basis of the orthogonal complement $W_j$ of $V_j$ in $V_{j+1}$, which satisfies
\begin{equation}
    V_j\oplus W_j = V_{j+1}, \forall j\in\mathbb{Z},
\end{equation}
where $\oplus$ means direct sum of two spaces.

We thus can decompose $f(t)$ in the orthonormal basis of $L^2(\mathbb{R})$
\begin{equation}
    f(t) = \sum_{k\in\mathbb{Z}}<f,\phi_{0,k}>\phi_{0,k}(t) + \sum_{j\geq 0}\sum_{k\in\mathbb{Z}}<f,\psi_{j,k}>\psi_{j,k}(t),
\end{equation}
where $<\cdot,\cdot>$ denotes inner product of two functions. $<f,\phi_{0,k}>$ and $<f,\psi_{j,k}>$ are wavelet coefficients.  
For the construction of orthonormal scaling function, one can refer to \cite{Daubechies1988}.

The multiresolution analysis described above for one input variable can be easily extended to the multi-input case. We present here a brief description of the extension to two-dimensions. The extension to higher dimensionality is straightforward. Defining the set of vector spaces $V^2_j\subset L^2(\mathbb{R}^2) $ by
\begin{equation}
     V^2_{j} = V_j\oplus V_j, \forall j\in\mathbb{Z},
\end{equation}
it is easy to show that the set $\{V^2_j\}_{j\in\mathbb{Z}}$ satisfies
\begin{equation*}
    V^2_j\subset V^2_{j+1},\quad
    \bigcap{V^2_j} = \{0\},\quad \overline{\bigcup V^2_j} = L^2(\mathbb{R}^2).
\end{equation*}

Defining
\begin{equation}
\begin{split}
    \Phi_{j,n_1,n_2}(x_1,x_2) &= 2^j\Phi(2^jx_1-n_1,2^jx_2-n_2) \\
    &= \phi_{j,n_1}(x_1)\phi_{j,n_2}(x_2),
\end{split}
\end{equation}
then we have
\begin{equation}
   V^2_j = linear\ span\{ \Phi_{j,n_1,n_2}, n_1,n_2\in\mathbb{Z})\}.
\end{equation}

Defining
\begin{equation}\label{2dPsi} 
\begin{split}
\Psi^1_{j,n_1,n_2}(x_1,x_2) &= 2^j\Psi(2^jx_1-n_1,2^jx_2-n_2) \\
&= \phi_{j,n_1}(x_1)\psi_{j,n_2}(x_2), \\
\Psi^2_{j,n_1,n_2}(x_1,x_2) &= 2^j\Psi(2^jx_1-n_1,2^jx_2-n_2) \\
&= \psi_{j,n_1}(x_1)\phi_{j,n_2}(x_2), \\
\Psi^3_{j,n_1,n_2}(x_1,x_2) &= 2^j\Psi(2^jx_1-n_1,2^jx_2-n_2) \\
&= \psi_{j,n_1}(x_1)\psi_{j,n_2}(x_2),
\end{split}
\end{equation}
and
\begin{equation}\label{2dWj} 
\begin{split}
W^{2,1}_j &= V_j\oplus W_j \\
&= linear\ span\{ \Psi^1_{j,n_1,n_2}, \ n_1,n_2\in\mathbb{Z})\}, \\
W^{2,2}_j &= W_j\oplus V_j \\
&= linear\ span\{ \Psi^2_{j,n_1,n_2}, \ n_1,n_2\in\mathbb{Z})\}, \\
W^{2,3}_j &= W_j\oplus W_j \\
&= linear\ span\{ \Psi^3_{j,n_1,n_2}, \ n_1,n_2\in\mathbb{Z})\},
\end{split}
\end{equation}
then we can obtain the orthogonal complement $W^2_j$ of $V^2_j$ in $V^2_{j+1}$, which satisfies
\begin{equation}
\begin{split}
    W^2_j &= W^{2,1}_j \oplus W^{2,2}_j \oplus W^{2,3}_j \\
    &= (V_j\oplus W_j)\oplus(W_j\oplus V_j)\oplus(W_j\oplus W_j),\\
    V^2_{j+1} &= V^2_j\oplus W^2_j, \forall j\in\mathbb{Z}.
\end{split}
\end{equation}
It is easy to show that the set $\{V^2_j,W^2_j\}_{j\in\mathbb{Z}}$ satisfies properties of the multiresolution analysis of $L^2(\mathbb{R}^2)$.

\subsection{Multiresolution Wavelet Neural Networks}
\begin{figure}[!t]
\centering
	\subfigure[]{
		\includegraphics[width=0.4\textwidth]{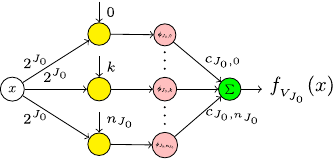}
	}
	\subfigure[]{
		\includegraphics[width=0.4\textwidth]{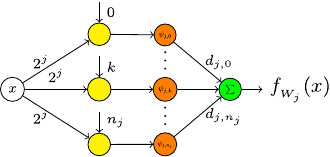}
	}
	
	\subfigure[]{
		\includegraphics[width=0.4\textwidth]{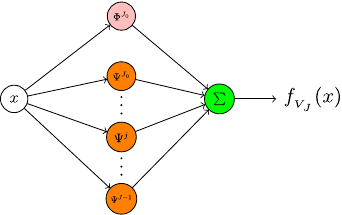}
	}
	\subfigure[]{
		\includegraphics[width=0.4\textwidth]{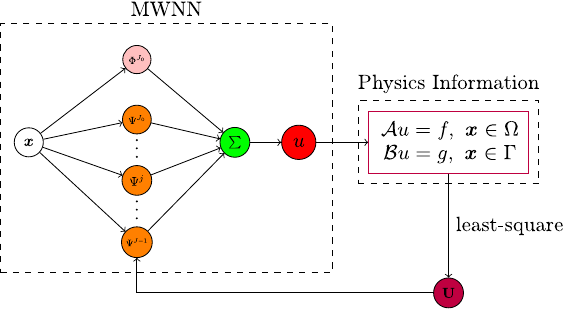}
	}
\caption{Structure of multiresolution wavelet neural networks and schematic of PIMWNN. In the subfigures (a) and (b) , the yellow nodes represent the scale transformation. The pink nodes represent the scaling function node and the orange nodes represent the wavelet function node. In the subfigure (c), the pink nodes represent the subnetwork of scaling function and the orange nodes represent subnetwork of wavelet function. The green nodes represent the summation operation. (a) subnetwork of scaling function for ${J_0}$, noted by $\Phi^{J_0}$. (b) subnetwork of wavelet function for ${j}$, noted by $\Psi^{j},j=J_0,\cdots,J-1$. (c) multiresolution wavelet neural network. The schematic diagram of PIMWNN is shown in subfigure (d).}
\label{fig:fig_mwnn}
\end{figure}

 The structure of a multiresolution wavelet neural networks is shown in Figure \ref{fig:fig_mwnn}. 
There are two types of nodes in the hidden layer: scaling function nodes and wavelet function nodes. The subnetwork of scaling function for ${J_0}$ constructed by scaling function nodes shown in Figure \ref{fig:fig_mwnn}(a) forms the orthonormal basis for the approximation of the function $f(x)$ at the ${J_0}$-th resolution. This subnetwork yields the following approximation
\begin{eqnarray}
    f(x)\approx f_{_{V_{J_0}}}(x) = \sum_{k=0}^{n_{{J_0}}}c_{_{J_0,k}}\phi_{{J_0,k}}(x).
\end{eqnarray}
Figure \ref{fig:fig_mwnn}(b) shows the subnetwork of wavelet function constructed by wavelet function nodes, which form the orthonormal basis for the approximation of the function $f(x)$ at the ${j}$-th resolution
\begin{eqnarray}
    f(x)\approx f_{_{W_{j}}}(x) = \sum_{k=0}^{n_j}d_{j,k}\psi_{j,k}(x).
\end{eqnarray}
By combining the scale function subnetwork and all wavelet function subnetworks, we obtain the multiresolution wavelet network shown in  Figure \ref{fig:fig_mwnn}(c), which forms the orthonormal basis for the approximation of the function $f(x)$ from the ${J_0}$-th to ${J}$-th resolutions
\begin{eqnarray}\label{eq_mra} 
\begin{split}
    f(x) &\approx f_{_{V_{J}}}(x) \\
    &= \sum_{k=0}^{n_{J_0}}c_{_{J_0,k}}\phi_{{J_0,k}}(x) + \sum^{J-1}_{j=J_0}\sum_{k=0}^{n_j}d_{j,k}\psi_{j,k}(x).
\end{split}
\end{eqnarray}

In this paper, we use the Shannon wavelet to construct the MWNN. Suppose the domain of function $f(x)$ is $[0,L]$, the Shannon scaling function and wavelet function are given by
\begin{subequations}
\begin{align}
\phi_{j,k}(x) &= 2^{j/2}\mathtt{sinc}\left(\frac{2^j}{L}x - k\right), \quad k = 0,1,\cdots,2^j, \\
\psi_{j,k}(x) &= 2^{j/2}\mathtt{sinc}\left(\frac{1}{2}(\frac{2^j}{L}x - k)\right)\cdot \cos\left(\frac{3\pi}{2}(\frac{2^j}{L}x - k)\right).
\end{align}
\end{subequations}
where $\mathtt{sinc}(x) = \sin(\pi x)/\pi x$. 

Assuming $J_0$ and $J$ are given, let $n_j = 2^j (j = J_0,\cdots,J-1)$, then the function $f(x)$ represented by the MWNN can be expressed as
\begin{eqnarray}\label{eq_mra_1}
    f(x)= \sum_{i=1}^N a_i\theta_i(x),
\end{eqnarray}
where, $a_i$ are the weights, $N=n_{J_0}+1+\sum_{j=J_0}^{J-1}(n_{j}+1)=2^{J}+J-J_0+1$, and $\theta_i$ represents the Shannon scaling or wavelet functions. 

\section{PIMWNN method}\label{method}
In this section, we introduce the basic concepts of PIMWNN method and demonstrate its applications to solve linear and nonlinear PDEs. The PIMWNN method integrates two ideas originating from PINN and MWNN. The first idea, inspired by PINN, is to encapsulate the physics of PDE within the loss function of the machine learning process. The other one coming from MWNN is to train the weights of the outer layer only. Since one of the difficulties that we found in the implementation of PINN was its slow training, we integrate PINN with MWNN to combine their advantages. Because we use the least-squares method to determine the weights of the outer layer, PIMWNN is faster than the iterative gradient descent-based methods.

\subsection{PIMWNN method for stationary problems}
Let us consider a PDE
\begin{subequations}
\begin{align}
\mathcal{A}u &= f,  \quad  {x}\in\Omega, \label{eq_lapace} \\
\mathcal{B}u &= g,  \quad  {x}\in\Gamma,    \label{eq_bou}
\end{align}
\end{subequations}
where $\Omega$ is a bounded domain in $\mathbb{R}^d$ with boundary 
$\Gamma=\Gamma_D\cup \Gamma_N$, and $\Gamma_D \cap \Gamma_N = \emptyset $. 
For example, if the differential operator $\mathcal{A}$ is given by
\begin{eqnarray}\label{eq_op_A}
    \mathcal{A}:= -\nabla\cdot (\alpha({{x}}) \nabla),
\end{eqnarray}
where $0<\alpha_0<\alpha({x})<\alpha_1$, \eqref{eq_lapace} is a diffusion equation. Here, $\alpha_0$ and $\alpha_1$ are constants. The boundary condition \eqref{eq_bou} can be Dirichlet, Neumann, or Robin types. In this paper, we consider Dirichlet boundary conditions: 
\begin{eqnarray}\label{eq_bouGD}
    u({x}) = g_D({x}), \quad x\in\Gamma_D.
\end{eqnarray}

According to equation \eqref{eq_mra_1}, we suppose the solution $u(x)$ can be expressed by the MWNN as
\begin{eqnarray}\label{eq_mra_u}
    u({x})= \sum_{i=1}^N u_i\theta_i(x)=\mathbf{\Theta}(x)\mathbf{{U}},
\end{eqnarray}
where $\mathbf{\Theta}(x) = (\theta_1(x),\theta_2(x),\cdots,\theta_N(x))$, $\mathbf U = (u_1,u_2,\cdots,u_N)^T$. Substituting \eqref{eq_mra_u} into \eqref{eq_lapace} and \eqref{eq_bouGD}, we can get
\begin{subequations}
\begin{align}
\mathcal{A}\mathbf{\Theta}(x)\mathbf{U} &= f(x), \qquad x\in\Omega, \\
\mathbf{\Theta}(x)\mathbf{U} &= g_D(x),        \,\,\quad  x\in\Gamma_D.
\end{align}
\end{subequations}
Taking some collocation points $\{x_i\}_{i=1}^{N_f}$ in $\Omega$ and $\{x_k\}_{k=1}^{N_b}$ on $\Gamma_D$, we get
\begin{eqnarray}\label{eq_lsq1}
    \begin{bmatrix}\mathbf{A} \\ \mathbf{B} \end{bmatrix} \mathbf{U} = \begin{bmatrix} \mathbf{F} \\\mathbf{G} \end{bmatrix},
\end{eqnarray}
where 
\begin{eqnarray*}
\begin{split}
&\mathbf{A} = (A_1,\cdots,A_i,\cdots,A_{N_f})^T, A_{i}=\mathcal{A}\mathbf{\Theta}(x_i),\\
&\mathbf{B} = (B_1,\cdots,B_k,\cdots,B_{N_b})^T, B_{k}=\mathbf{\Theta}(x_k), \\
&\mathbf{F} = (f(x_1),\cdots,f(x_i),\cdots,f(x_{N_f}))^T, \\
&\mathbf{G} = (g_D(x_1),\cdots,g_D(x_k),\cdots,g_D(x_{N_b}))^T. 
\end{split}
\end{eqnarray*}
Then we can obtain the weights vector $\mathbf{U}$ by solving a least-squares problem with the linear system \eqref{eq_lsq1}.

The algorithm for the PIMWNN method is summarized as Algorithm 1. The schematic diagram of PIMWNN is shown in Fig. \ref{fig:fig_mwnn}(d) .

\begin{table}[h!]\label{stepwnn}
\begin{minipage}{15cm}
 \begin{tabular*}{11.5cm}{lc}
 \toprule
 \textbf{Algorithm 1:} Physics Informed Multiresolution Wavelet Neural Networks Method & \\ 
 \midrule
  \quad \textbf{Step 1}: Choose the coarsest resolution $J_0$ and the finest resolution $J$. & \\
  \quad \textbf{Step 2}: Construct the multiresolution wavelet neural network with $J_0$ and $J$. & \\
  \quad \textbf{Step 3}: Specify the train data set $\{\mathbf{x}_i\}^{N_f}_{i=1}$ in $\Omega$ and $\{\mathbf{x}_k\}_{k=1}^{N_b}$ on $\Gamma_D$. &\\
  \quad \textbf{Step 4}: Assemble the linear system $\mathbf{AU=F}$ and $\mathbf{BU=G}$ in Eq. \eqref{eq_lsq1}.& \\
  \quad \textbf{Step 5}: solving the least-squares problem with the linear system \eqref{eq_lsq1} to obtain $\mathbf{U}$. \\
  \quad \textbf{Step 6}: Update the weights of the outer layer of the network with $\mathbf{U}$.\\
 \bottomrule
 \end{tabular*}
\end{minipage}
\end{table}

\subsection{PIMWNN method for time-dependent problems}\label{sc-time-dependent}
In this section, we use the PIMWNN method to solve time-dependent PDEs using
a space-time approach. We apply the PIMWNN method to solve time-dependent problems in a unified manner rather than discretizing time and space separately.

Let us consider the following time-dependent problem
\begin{subequations}
\begin{align}
\mathcal{A}u &= f,  \quad x\in\Omega\times I,      \label{eq_lapace_t} \\
\mathcal{B}u &= g,  \quad x\in\Gamma\times I,      \label{eq_bou_t} \\
\mathcal{C}u &= h,  \quad x\in\Omega\times \{0\},  \label{eq_ini_t} 
\end{align}
\end{subequations}
where $\Omega$ is a bounded domain in $\mathbb{R}^d$ with boundary $\Gamma=\Gamma_D\cup \Gamma_N$, and $\Gamma_D\cap \Gamma_N = \emptyset$, and $I=(0,T]$ is the time interval of interest. For example, if the differential operator $\mathcal{A}$ is given by
\begin{eqnarray}\label{eq_op_A_t}
    \mathcal{A}:= \frac{\partial}{\partial t} -\nabla\cdot (\alpha({x}) \nabla),
\end{eqnarray}
where $0<\alpha_0<\alpha(x)<\alpha_1$, \eqref{eq_lapace_t} is a heat conduction equation. Here, $\alpha_0$ and $\alpha_1$ are constants. The boundary condition \eqref{eq_bou_t} can be Dirichlet, Neumann, or Robin types. Here, we consider Dirichlet boundary conditions: 
\begin{eqnarray}\label{eq_bouGD_t}
    u(x,t) = g_D(x,t), \quad x\in\Gamma_D\times I,
\end{eqnarray}
and the initial condition \eqref{eq_ini_t} is given by
\begin{eqnarray}\label{eq_ini_t0}
    u(x,0) = h_0(x), \quad x\in\Omega\times \{0\}.
\end{eqnarray}

Unlike the traditional approach, for example, the finite difference method, discretizing temporal and
spatial variables separately, here the space-time approach treats the temporal and spatial variables equally and jointly.

We express the solution $u(x,t)$ by the MWNN as
\begin{equation}\label{eq_mra_u_t}
    u(x,t)= \sum_{i=1}^N u_i\theta_i(x,t)=\mathbf{\Theta}(x,t)\mathbf{{U}},
\end{equation}
where $\mathbf{\Theta}(x,t) = (\theta_1(x,t),\theta_2(x,t),\cdots,\theta_N(x,t))$, $\mathbf U = (u_1,u_2,\cdots,u_N)^T$. Substituting \eqref{eq_mra_u_t} into \eqref{eq_lapace_t} and \eqref{eq_bouGD_t}, we can get
\begin{subequations}
\begin{align}
\mathcal{A}\mathbf{\Theta}(x,t)\mathbf{U} &= f(x,t), \quad x\in\Omega\times I, \\
\mathbf{\Theta}(x,t)\mathbf{U} &= g_D(x,t),           \quad x\in\Gamma_D\times I, \\
\mathbf{\Theta}(x,0)\mathbf{U} &= h_0(x),       \quad x\in\Omega,
\end{align}
\end{subequations}
taking some collocation points $\{(x_i,t_i)\}_{i=1}^{N_f}$ in $\Omega\times I$, $\{(x_k,t_k)\}_{k=1}^{N_b}$ on $\Gamma_D\times I$ and $\{(x_m,0)\}_{m=1}^{N_t}$ in $\Omega$, we get
\begin{equation}\label{eq_lsq2}
    \begin{bmatrix}\mathbf{A} \\ \mathbf{B} \\ \mathbf{C}\end{bmatrix} \mathbf{U} =
    \begin{bmatrix} \mathbf{F} \\ \mathbf{G} \\ \mathbf{H}\end{bmatrix},
\end{equation}
where 
\begin{equation}
\begin{split}
\mathbf{A} &= (A_1,\cdots,A_i,\cdots,A_{N_f})^T,  A_{i}=\mathcal{A}\mathbf{\Theta}(x_i,t_i),\\
\mathbf{B} &= (B_1,\cdots,B_k,\cdots,B_{N_b})^T, B_{k}=\mathbf{\Theta}(x_k,t_k),\\
\mathbf{C} &= (C_1,\cdots,C_m,\cdots,C_{N_t})^T, C_{m}=\mathbf{\Theta}(x_m,0),\\
\mathbf{F} &= (f_1,\cdots,f_i,\cdots,f_{N_f})^T,     f_i = f(x_i,t_i),\\
\mathbf{G} &= (g_1,\cdots,g_k,\cdots,g_{N_b})^T, g_k = g_D(x_k,t_k),\\
\mathbf{H} &= (h_1,\cdots,h_m,\cdots,h_{N_t})^T, h_m = h_0(x_m).
\end{split}
\end{equation}
Finally, we can obtain the weights vector $\mathbf{U}$ by solving a least-squares problem with the linear system \eqref{eq_lsq2}.

\subsection{PIMWNN method for time-dependent nonlinear problems}\label{sc-bergers}
In previous sections, we demonstrated the use of PIMWNN method to solve linear PDEs. Next, we use Burger's equation as an example to show how to apply PIMWNN method to solve nonlinear PDEs. The 1-D Burger's equation is given as follows:
\begin{subequations}\label{eq_bur}
\begin{align}
&u_t + uu_x -\epsilon u_{xx} = f, \quad  \quad  x\in\Omega\times I,      \label{eq_bur_t} \\
&u(x,t) = g_D(x,t), \quad \qquad  \ x\in\Gamma_D\times I,      \label{eq_bou_bur} \\
&u(x,0) = h_0(x),   \ \ \, \quad \qquad  \ x\in\Omega\times \{0\}.  \label{eq_ini_bur} 
\end{align}
\end{subequations}
To address the above problem, we can apply nonlinear iteration methods such as Picard iteration \cite{Langtangen2017SolvingNO}. Using backward Euler scheme for \eqref{eq_bur_t}, we can get
\begin{equation}
    \frac{u^n - u^{n-1}}{dt} = -u^n u_x^n +\epsilon u_{xx}^n + f.
\end{equation}
We introduce a Picard iteration with k as iteration counter. Let $U^{k} = u^{n,k}$. A linearization of the $u^nu_x^n$ term in iteration $k+1$ is to use the previously computed $U^{k}$ approximation in the convection coefficient: $u^n=U^{k}$. The unknown function $U^{k+1}$ then fulfills the linear PDE
\begin{equation}
    \frac{U^{k+1} - u^{n-1}}{dt} = -U^{k} U_x^{k+1} +\epsilon U_{xx}^{k+1} + f.
\end{equation}
The initial guess for the Picard iteration at this time level can be taken as the solution at the
previous time level: $U^{0} = u^{n-1}$. When k reaches the maximum number of iteration steps $K$, let $u^{n} = U^{K}$.

\section{Numerical examples}\label{num_exp}
This section presents several examples to demonstrate the performance of PIMWNN method. Firstly, we test the representation power of PIMWNN by a 1D composite function and a 2D Gaussian function. Secondly, we rigorously test PIMWNN method on various linear, quasi-linear and non-linear PDEs. All these examples are taken from \cite{Berg2017AUD,Dwivedi2019PhysicsIE,raissi2019}. We evaluate the accuracy of numerical solution by relative $L^2$ error. All the experiments are conducted in python environment running in an Intel Core i7 2.20 GHz CPU and 16 GB RAM Thinkpad laptop.

\subsection{Representation of functions with sharp gradients}
\begin{figure}[!t]
\centering
	\subfigure[Solution ($J$ = 11)]{
		\includegraphics[width=0.36\textwidth]{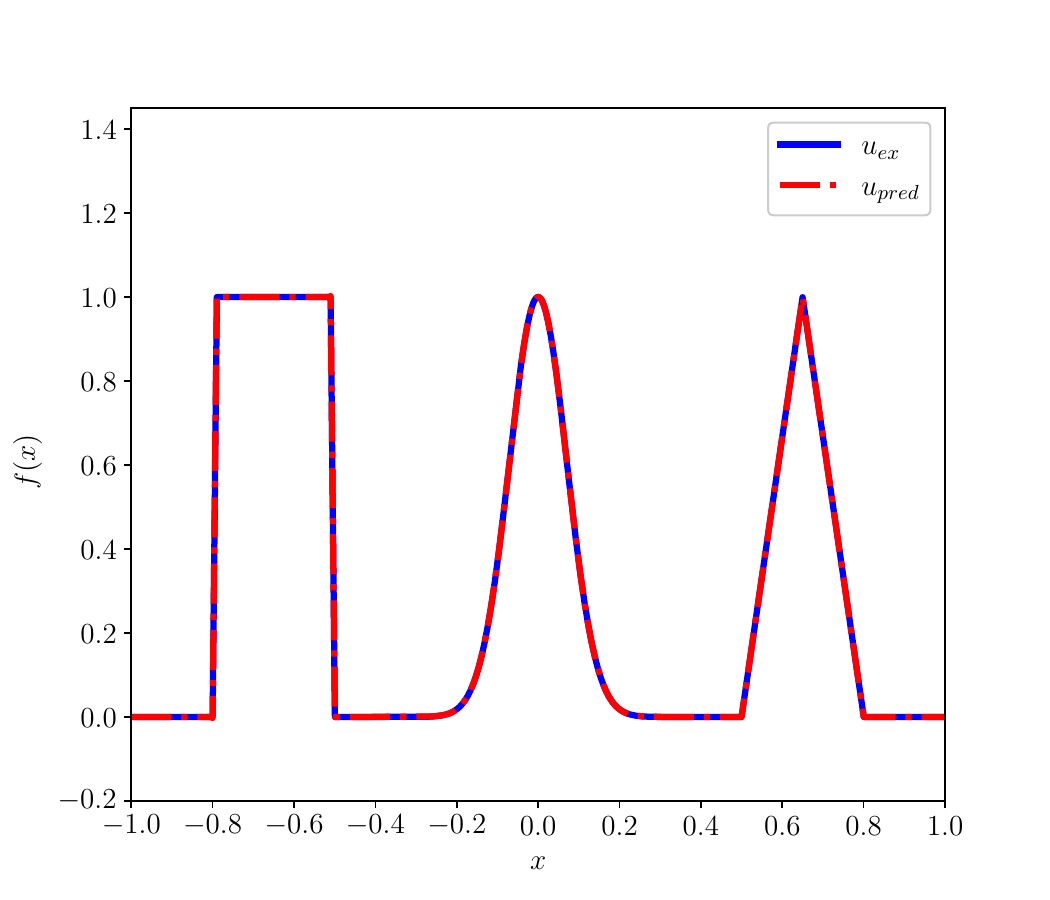}
	}
	\subfigure[$J$ = 9]{
		\includegraphics[width=0.4\textwidth]{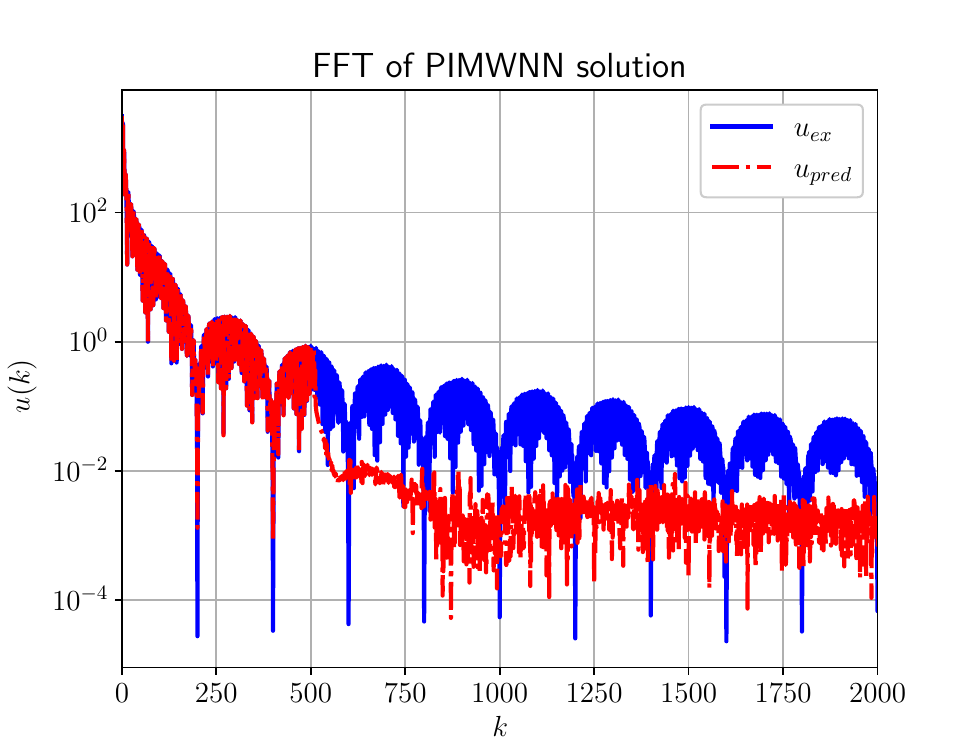}
	}
	
	\subfigure[$J$ = 10]{
		\includegraphics[width=0.4\textwidth]{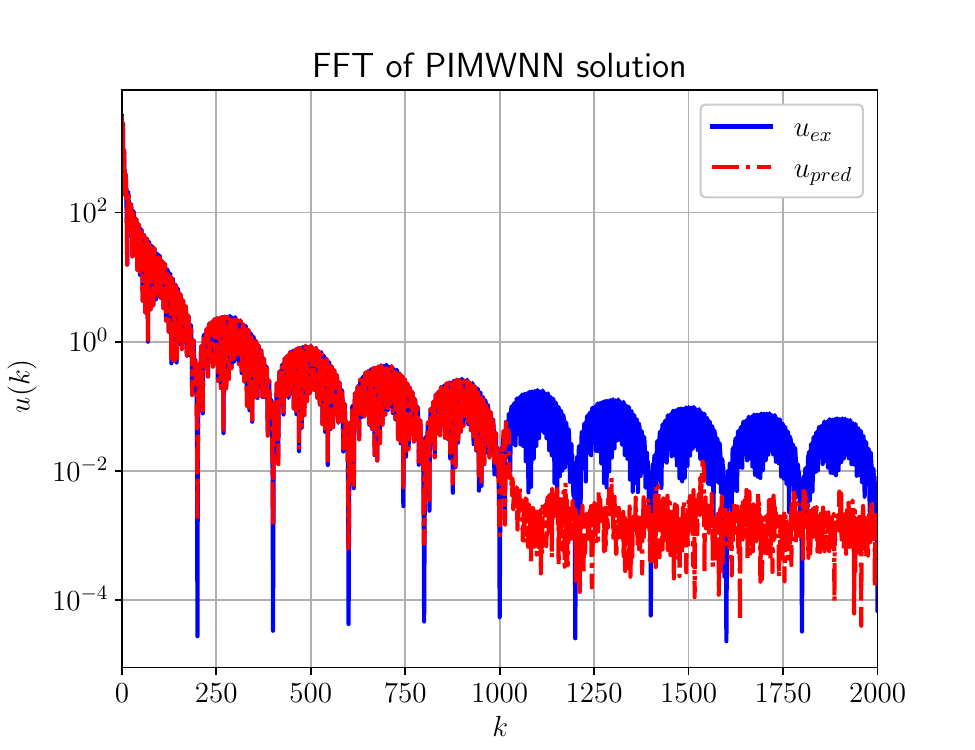}
	}
	\subfigure[$J$ = 11]{
		\includegraphics[width=0.4\textwidth]{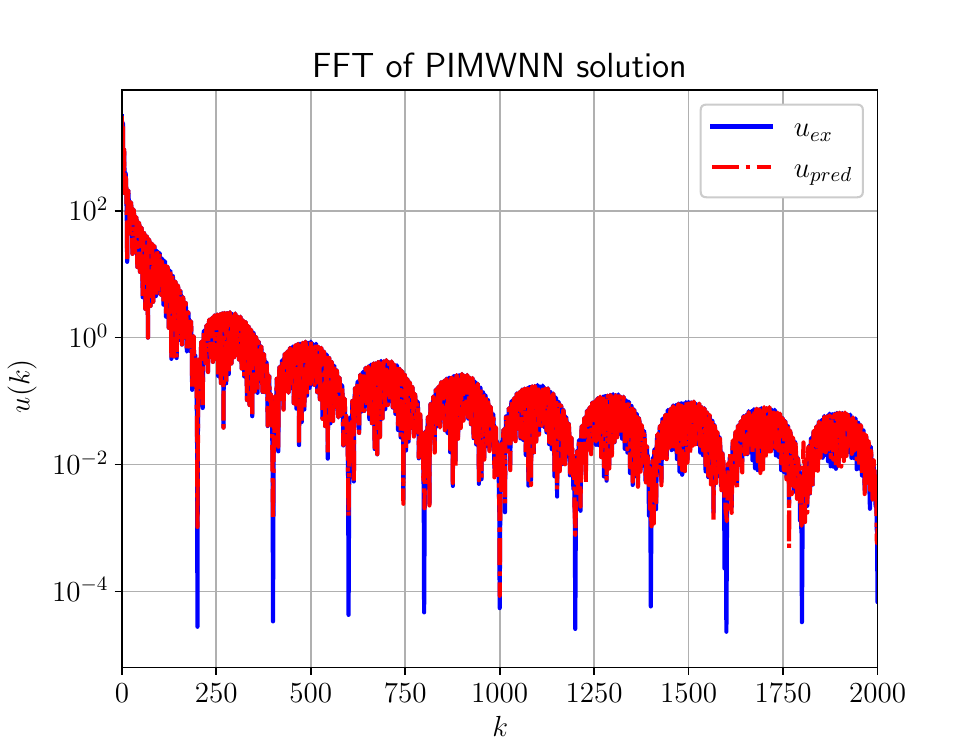}
	}
\caption{Representation of function $f_1(x)$ and FFT of PIMWNN solutions.}
\label{tc-0-1}
\end{figure}

The two functions with sharp gradients are given as
\begin{equation}
f_1(x) = 
\begin{cases}
\frac{1}{2}\{ sgn(x+0.8)-sgn(x+0.5) \} & x \in [-1,-\frac{1}{2}]  \\
\exp(-100x^2) & x \in (-\frac{1}{2},\frac{1}{2}]  \\
\frac{20}{3}x - \frac{10}{3} & x \in (\frac{1}{2},\frac{13}{20}]  \\
-\frac{20}{3}x + \frac{16}{3} & x \in (\frac{13}{20},\frac{16}{20}]  \\
0 & x \in (\frac{16}{20},1],
\end{cases}
\end{equation}
and
\begin{equation}
f_2(x,y) = \exp(-20(x^2+y^2)).
\end{equation}
The two functions were employed to evaluate the representation capability of PIELM \cite{Dwivedi2019PhysicsIE}. As illustrated in Figures 17 and 18 of reference \cite{Dwivedi2019PhysicsIE}, PIELM exhibits limitations in accurately representing profiles with sharp gradients and corners.. Here, we use these two functions to test the representation power of PIMWNN method also. The function $f_1(x)$ and its  PIMWNN solution are shown in Fig. \ref{tc-0-1}(a). In contrast to the noticeable discrepancy between the exact values and PIELM predictions observed in \cite{Dwivedi2019PhysicsIE}, our PIMWNN predictions demonstrate excellent agreement with the exact values. 
The results unequivocally demonstrate that the PIMWNN method excels in representing functions with sharp gradients and corners. 

In PINNs, the spectral bias inherent in neural networks poses a significant challenge. This issue is a well-documented phenomenon where neural networks exhibit a tendency to learn lower frequencies significantly faster than higher frequencies \cite{Xuz2020,Rahaman2019}. To analyze the spectral characteristics of the solutions, we conducted fast Fourier transforms (FFT) on both the exact solution and the PIMWNN solution depicted in Fig. \ref{tc-0-1}(a), and the results are presented in Figs. \ref{tc-0-1}(b)-(d). In these figures, $u(k)$ and $k$ are the amplitude and wave number in arbitrary units, respectively. The FFTs of the PIMWNN solution of the function $f_1(x)$ with scale $J$ from 9 to 11 are given in Figs. \ref{tc-0-1}(b)-(d). It can be seen from Fig. \ref{tc-0-1}(b) that the spectra of PIMWNN solution with the scale $J$=9 and the exact solution can match very well with each other before $k$=500. To improve the high-frequency accuracy of the PIMWNN solution, we incrementally increased the scale $J$, and find that spectra of PIMWNN solution and exact solution can match very well with each other before $k$=1000 with $J$=10 (Fig. \ref{tc-0-1}(c)) and before $k$=2000 with $J$=11 (Fig. \ref{tc-0-1}(d)). Consequently, the PIMWNN solution exhibits remarkable congruence with the exact solution.

\begin{figure}[!t]
\centering
	\subfigure[function $f_2(x,y)$]{
		\includegraphics[width=0.3125\textwidth]{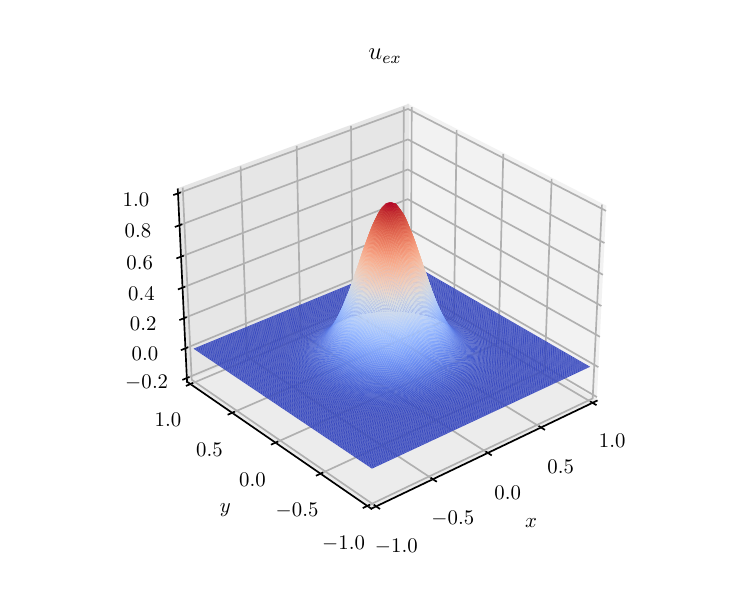}
	}
	\subfigure[PIMWNN solution]{
		\includegraphics[width=0.3125\textwidth]{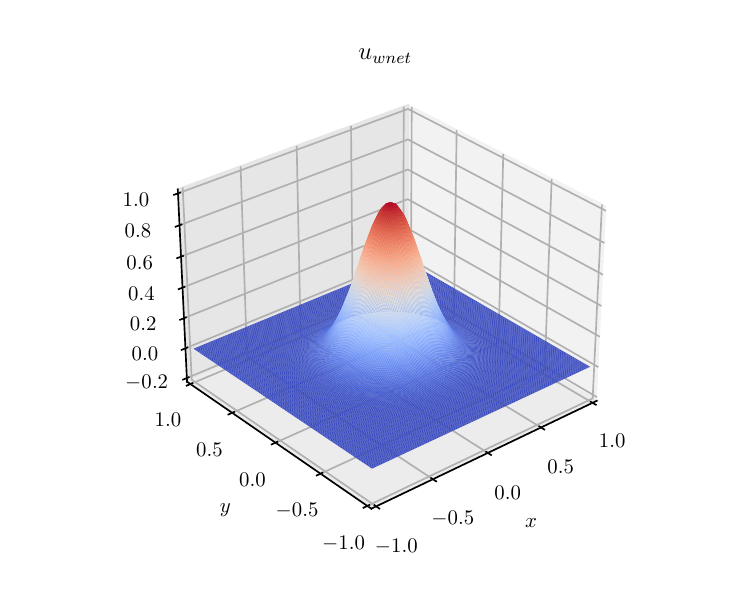}
	}	
	\subfigure[Pointwise error]{
		\includegraphics[width=0.25\textwidth]{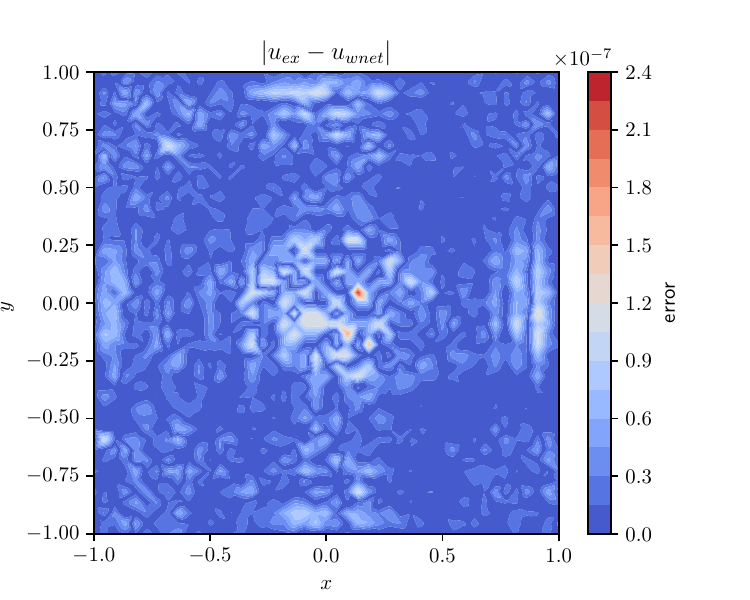}
	}
\caption{Representation of function $f_2(x,y)$.}
\label{tc-0-2}
\end{figure}

The function $f_2(x,y)$ is illustrated in Fig. \ref{tc-0-2}(a), the PIMWNN solution of  $f_2(x,y)$ with $J^x\times J^y=3\times 3$ is depicted in Fig. \ref{tc-0-2}(b). The pointwise error between the function $f_2(x,y)$ and the PIMWNN solution is shown in Fig. \ref{tc-0-2}(c), where we can see that the relative $L^2$ error is on order of magnitude $\mathcal{O}(10^{-7})$. Similarly, FFT analysis can be conducted for the function $f_2(x,y)$. 

These two examples are summarized in Tables \ref{btc-0-1} and \ref{btc-0-2}, we can observe that PIMWNN method provides more accurate numerical solutions with scale $J$ from 7 to 11 for the function $f_1(x)$ and from $1\times 1$ to $3\times 3$ for the function $f_2(x,y)$, given the same data of other parameters. 

\subsection{1D steady cases}
\subsubsection{1D advection, diffusion and dvection-diffusion equation}
\begin{figure}[!t]
\centering
	\subfigure[Solution of Eq. (\ref{Tc-1}).]{
		\includegraphics[width=0.35\textwidth]{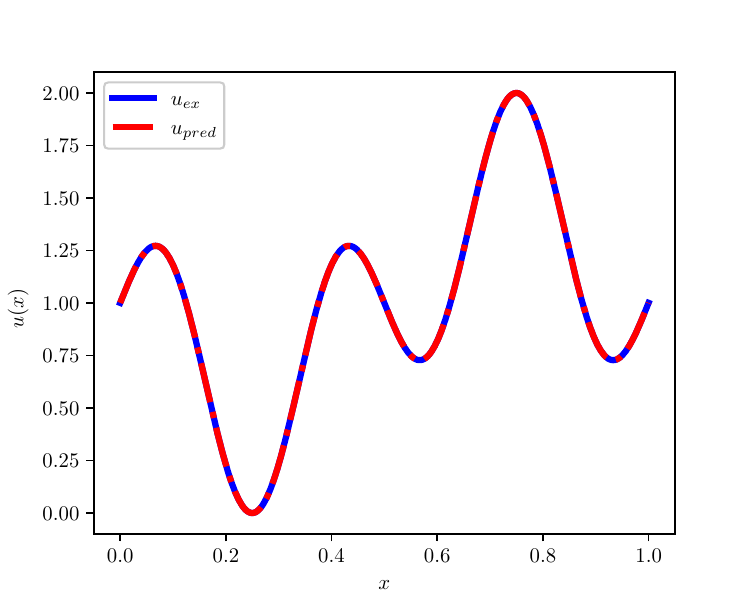}
	}
	\subfigure[Solution pointwise error of Eq. (\ref{Tc-1}).]{
		\includegraphics[width=0.35\textwidth]{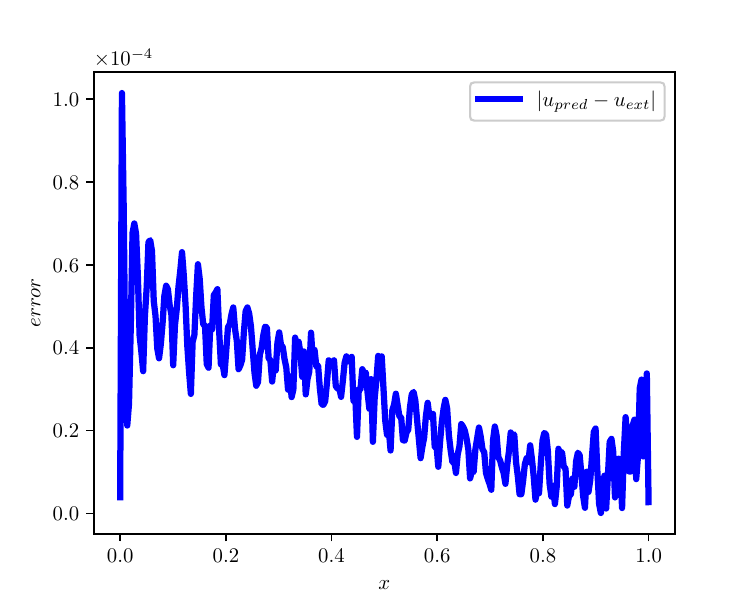}
	}	
	
	\subfigure[Solution of Eq. (\ref{Tc-2}).]{
		\includegraphics[width=0.35\textwidth]{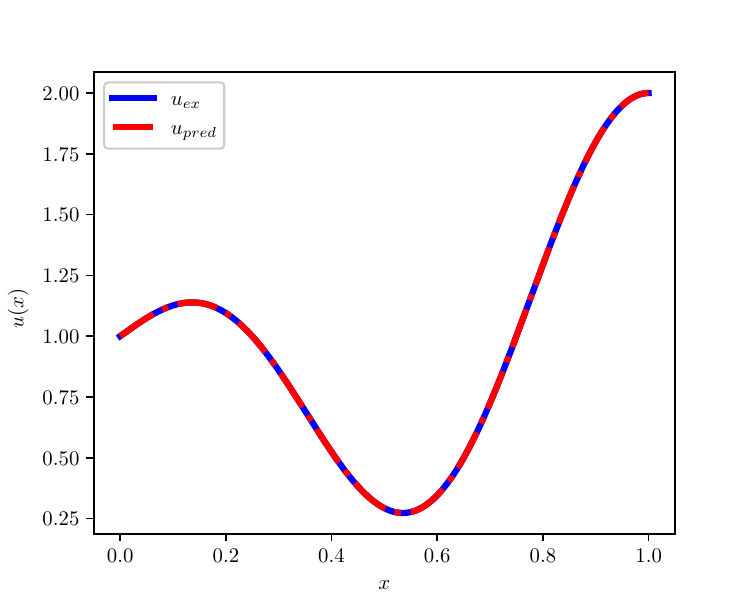}
	}
	\subfigure[Solution pointwise error of Eq. (\ref{Tc-2}).]{
		\includegraphics[width=0.35\textwidth]{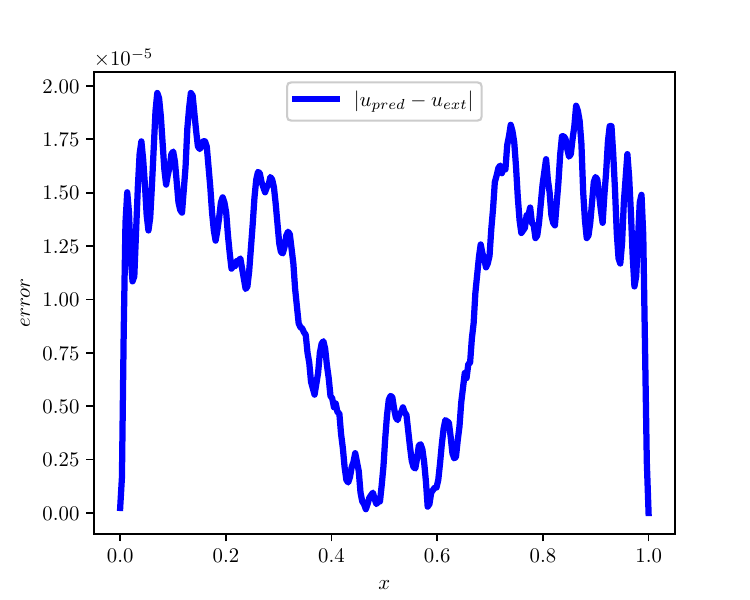}
	}	
	
	\subfigure[Solution of Eq. (\ref{Tc-3}).]{
		\includegraphics[width=0.35\textwidth]{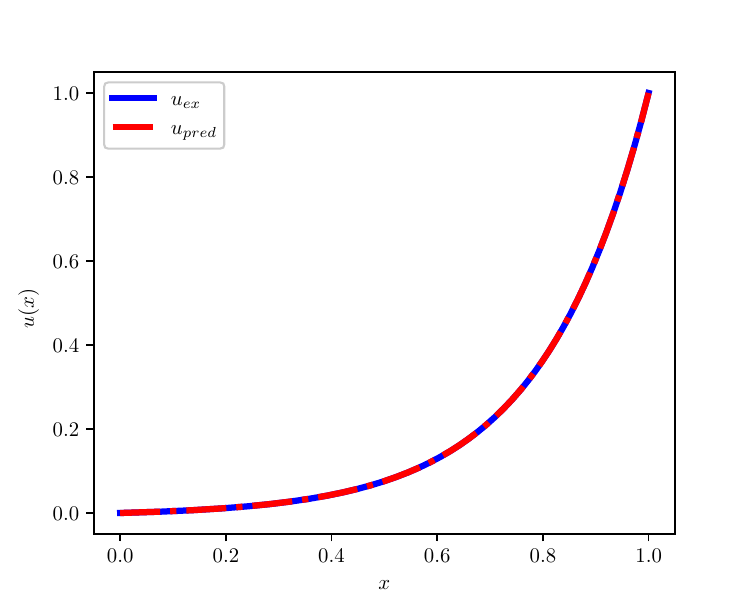}
	}
	\subfigure[Solution pointwise error of Eq. (\ref{Tc-3}).]{
		\includegraphics[width=0.35\textwidth]{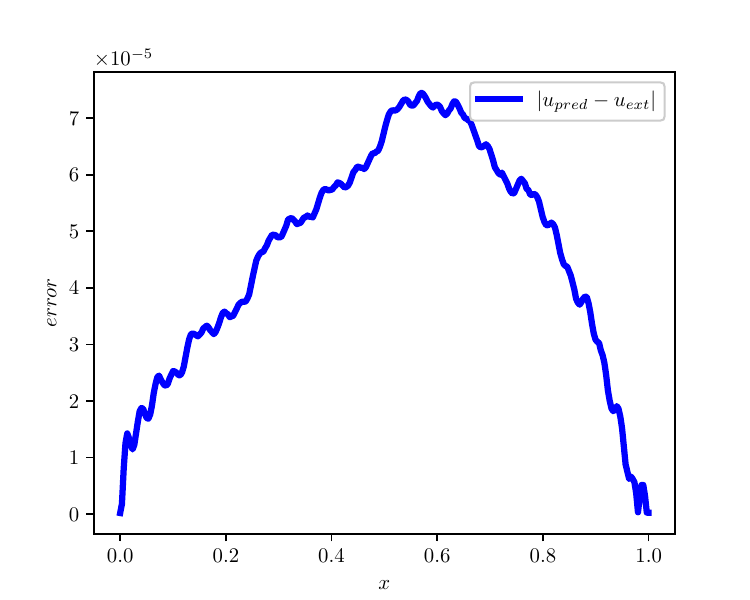}
	}	
\caption{Solution and error plots of equations (\ref{Tc-1})-(\ref{Tc-3}).}
\label{tc-1-3}
\end{figure}

In these examples, we solve 1D stationary advection equation (\ref{Tc-1}), diffusion equation (\ref{Tc-2}) and advection-diffusion equation (\ref{Tc-3}), given by
\begin{subequations}
\begin{align}
&u_x = f, \quad 0 < x < 1,   \label{Tc-1}  \\
&u_{xx} = f, \quad 0 < x < 1,   \label{Tc-2} \\
&u_x - vu_{xx} = f, \quad 0 < x < 1.   \label{Tc-3} 
\end{align}
\end{subequations}
The Dirichlet boundary conditions and the function $f$ are chosen such that the exact solution $u$ are given by
\begin{subequations}
\begin{align}
u(x) &= \sin(2\pi x)\cos(4\pi x) + 1, \\
u(x) &= \sin(\frac{\pi x}{2})\cos(2\pi x) + 1, \\
u(x) &= \frac{\exp{\frac{x}{v}}-1}{\exp{\frac{1}{v}}-1},
\end{align}
\end{subequations}
respectively. The results for these test cases are given in Fig. \ref{tc-1-3}. The left figure shows the exact solution and PIMWNN solution. The right figure shows the absolute value of their differences. The summary of the experiments of Eq. (\ref{Tc-1}), Eq.(\ref{Tc-2}) and Eq.(\ref{Tc-3}) are given in Table \ref{Tc1-3}. We can still observe that PIMWNN method provides more accurate numerical solutions with a larger number of scale $J$, given the same data of other parameters. 

\begin{table}[!t]
\caption{Summary of experiments for function $f_1(x)$.}\label{btc-0-1}
\centering
\begin{tabular}{cccccccc}
\hline
  Fucntion & $N_f$ & $N_b$ & $J_0$ & $J$ & $N$ & $e_{L^2}$ & $t_{cpu}(s)$ \\
\hline
          & 5000 & 2 & 0 &  7  &  265 & 1.055e-02 &  2.802e-02\\
  $f_1(x)$  & 5000 & 2 & 0 &  9  & 1035 & 1.426e-03 &  1.967e-01\\
          & 5000 & 2 & 0 & 11  & 4109 & 2.018e-04 &  1.053e+01\\
\hline
\end{tabular}
\end{table}

\begin{table*}[!t]
\caption{Summary of experiments for function $f_2(x,y)$.}\label{btc-0-2}
\centering
\begin{tabular}{cccccccc}
\hline
  Function  & $N_f$ & $N_b$ & $J^x_0\times J^y_0$ & $J^x\times J^y$ & $N$ & $e_{L^2}$ & $t_{cpu}(s)$ \\
\hline
          & 2000 & 400 & 0$\times$0 & 1$\times$1 & 7$\times$7 & 1.185e-01   & 4.343e-03\\
  $f_2(x,y)$  & 2000 & 400 & 0$\times$0 & 2$\times$2 & 12$\times$12 & 1.304e-03 & 1.159e-02\\
          & 2000 & 400 & 0$\times$0 & 3$\times$3 & 21$\times$21 & 1.181e-07 & 4.966e-02\\
\hline
\end{tabular}
\end{table*}

\begin{table}[!t]
\caption{Summary of experiments for equations (\ref{Tc-1})-(\ref{Tc-3}).}\label{Tc1-3}
\centering
\begin{tabular}{cccccccc}
\hline
  Equations & $N_f$ & $N_b$ & $J_0$ & $J$ & $N$ & $e_{L^2}$ & $t_{cpu}(s)$ \\
\hline
             & 100  & 2 & 0 & 1  &  7  & 3.175e-01 & 1.448e-02\\
Eq.(\ref{Tc-1}) & 100  & 2 & 0 & 2  & 12  & 7.927e-03 & 2.908e-02\\
             & 100  & 2 & 0 & 3  & 21  & 1.318e-03 & 6.015e-02\\
\hline
       & 100  & 2 & 0 & 1  &  7  & 9.821e-03 & 1.278e-02\\
  Eq.(\ref{Tc-2}) & 100  & 2 & 0 & 2  & 12  & 1.316e-03 & 2.748e-02\\
       & 100  & 2 & 0 & 3  & 21  & 1.688e-04 & 6.415e-02\\
\hline
       & 100  & 2 & 0 & 1  &  7  & 1.632e-01 & 1.248e-02\\
  Eq.(\ref{Tc-3}) & 100  & 2 & 0 & 3  & 21  & 2.709e-03 & 5.981e-02\\
       & 100  & 2 & 0 & 5  & 71  & 1.228e-04 & 2.881e-01\\
\hline
       & 20000  & 2 & 0 & 5  &  71  & 4.592e00 & 2.010e00\\
  Eq.(\ref{Tc-1d-helmholtz}) & 20000  & 2 & 0 & 6  & 136  & 2.356e-02 & 1.021e01\\
       & 20000  & 2 & 0 & 7  & 265  & 2.174e-03 & 8.031e01\\
\hline
\end{tabular}
\end{table}

\subsubsection{1D Helmholtz Equation}
\begin{figure}[!bt]
\centering
	\subfigure[Solution.]{
		\includegraphics[width=0.3\textwidth]{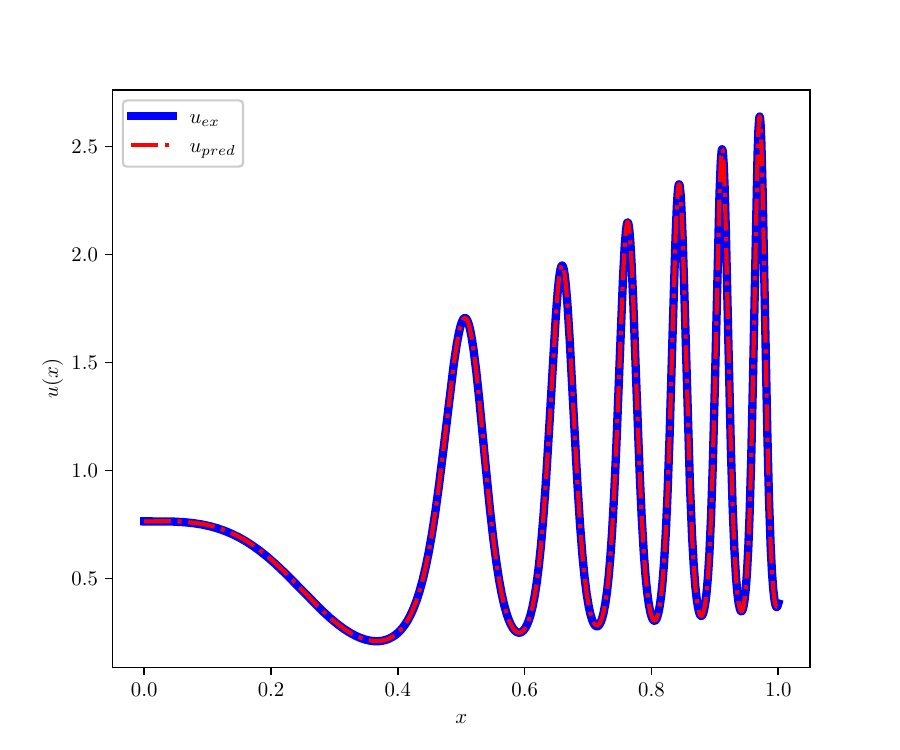}
	}
	\subfigure[Pointwise error.]{
		\includegraphics[width=0.3\textwidth]{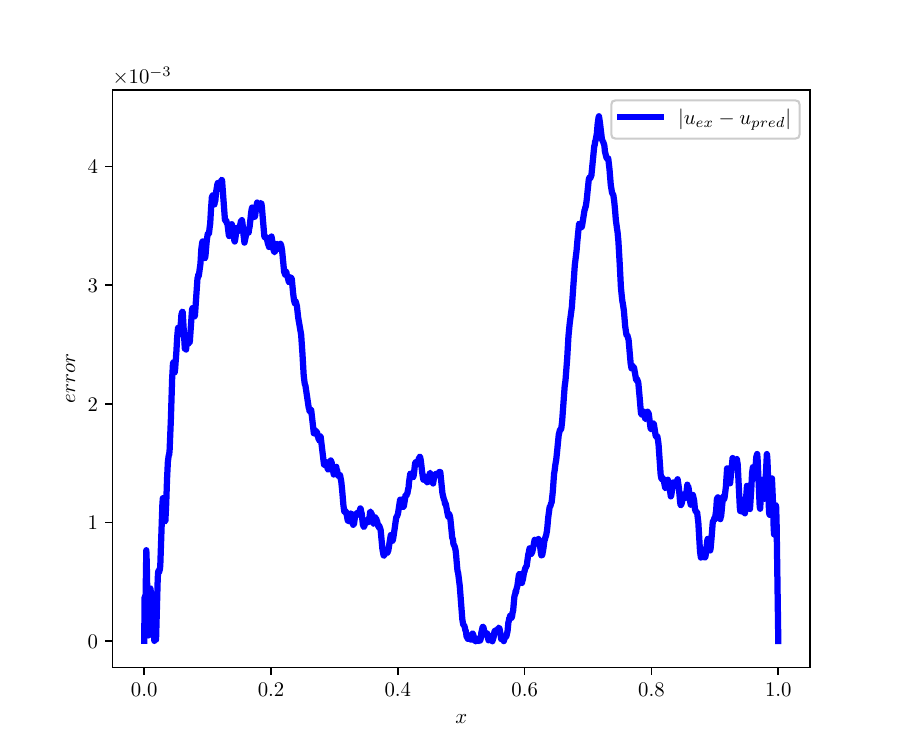}
	}	
	\subfigure[FFT of solution.]{
		\includegraphics[width=0.3\textwidth]{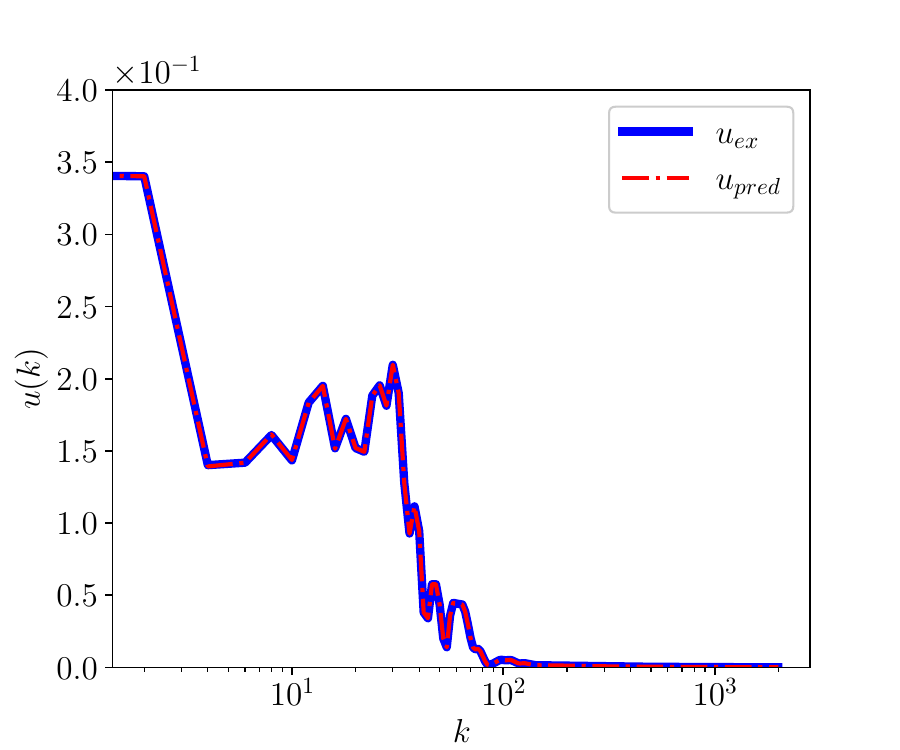}
	}	
\caption{Solution and error plots for 1D Helmholtz equation  Eq.(\ref{Tc-1d-helmholtz}).}\label{helmholtz1d}
\end{figure}

\begin{figure}[!bt]
\centering
	\subfigure[Solution.]{
		\includegraphics[width=0.28\textwidth]{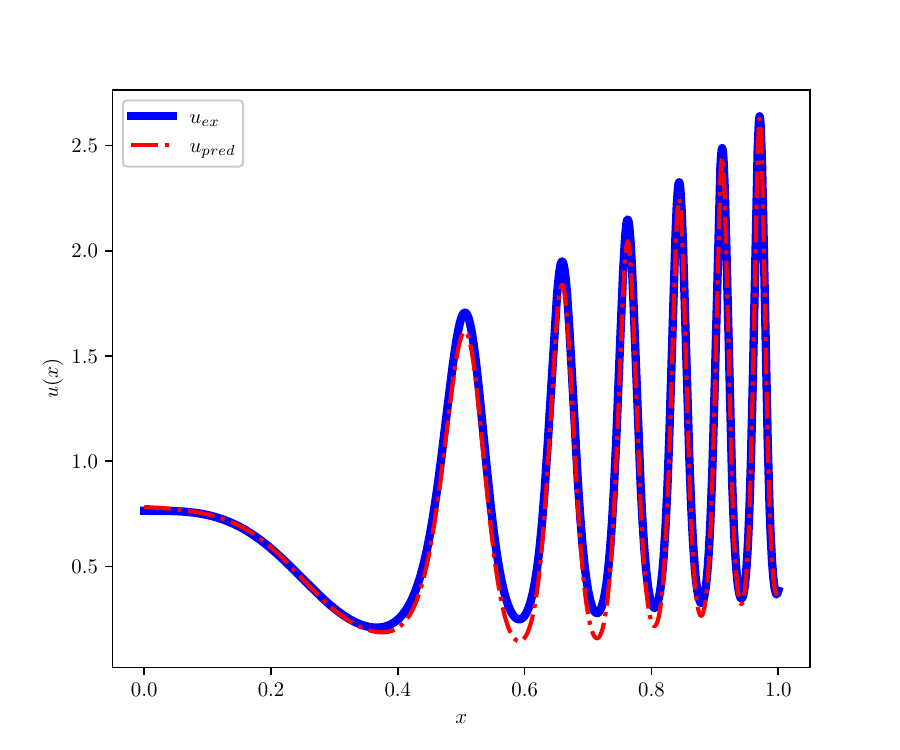}
	}
	\subfigure[Pointwise error.]{
		\includegraphics[width=0.31\textwidth]{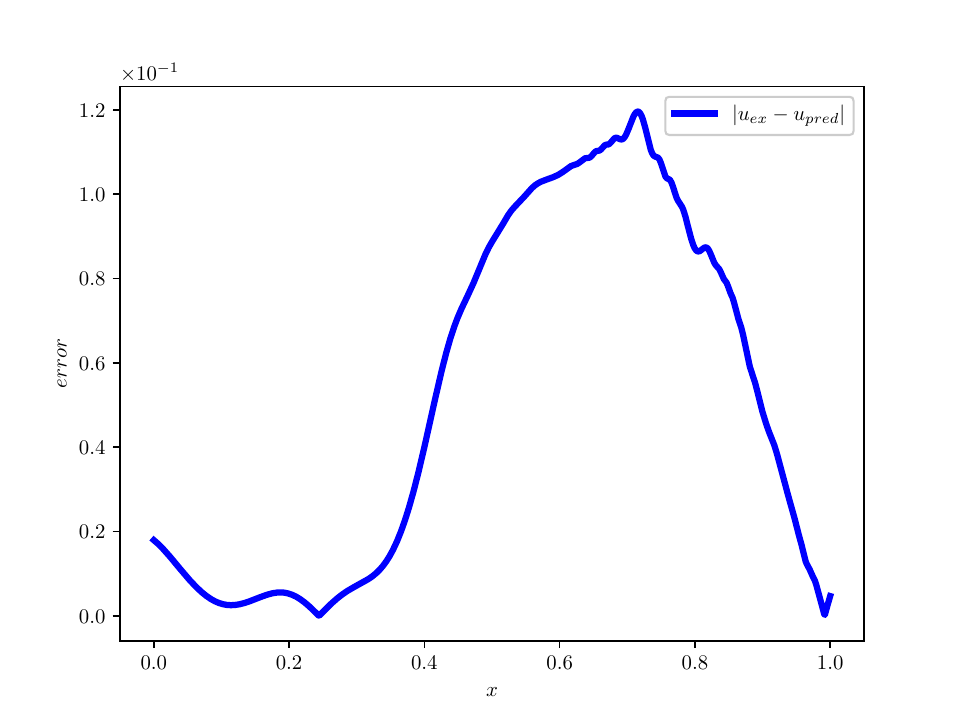}
	}	
	\subfigure[FFT of solution.]{
		\includegraphics[width=0.31\textwidth]{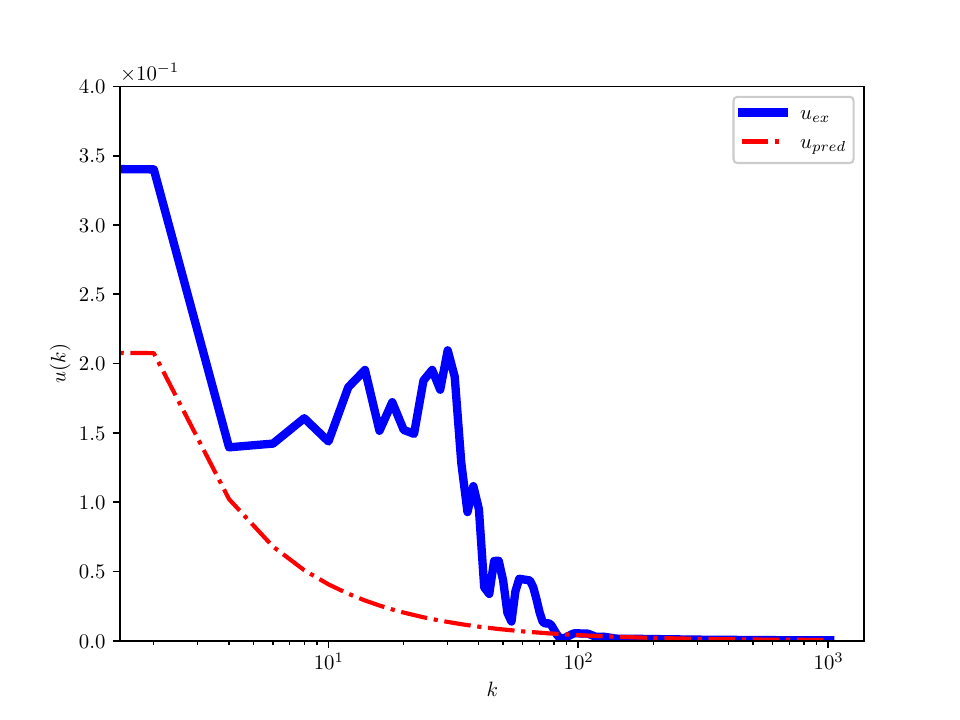}
	}	
\caption{PINN solution and error plots for 1D Helmholtz equation  Eq.(\ref{Tc-1d-helmholtz}).}\label{helmholtz1d-pinn}
\end{figure}

In this example, we solve a Helmholtz equation with high-frequency oscillations on the domain $\Omega=[0,1]$,
\begin{subequations}
\begin{align}
-u_{xx} + \lambda u &= f(x), \quad x \in (0,1),   \label{Tc-1d-helmholtz} \\
u(0) &= g_1(x), \\
u(1) &= g_2(x), 
\end{align}
\end{subequations}
where the $\lambda = 10$ and $f(x)$ is a source term, $g_1(x)$ and $g_2(x)$ are boundary conditions. $f(x), g_1(x)$ and $g_2(x)$ are chosen such that the exact solution are given by 
\begin{equation}
u(x) = \frac{1}{2}(x^2 + 1)e^{cos(40x^3-24)}.
\end{equation} 
The exact solution and the PIMWNN solution are shown in Fig. \ref{helmholtz1d}(a) and the pointwise error is shown in Fig. \ref{helmholtz1d}(b), where the scale $J = 7$. We can see that the the pointwise error is on order of magnitude $\mathcal{O}(10^{-3})$. the relative  error $L^2 = 6.887\times 10^{-2}$. The FFT of solution are shown in Fig. \ref{helmholtz1d}(c). The results show that the spectrum of the PIMWNN solution matches well with that of the exact solution. Solving the problem by the PIMWNN method took 80.31 seconds (see Table 3). For comparison, we use the PINN method to solve the above problem, and the PINN method is implemented with Package DeepXDE\cite{Lu2021deepxde}. The neural network of PINN consists of two hidden layers, each with 80 nodes. The activation function is the tanh function. The Adam algorithm is used for training, with a learning rate of 0.001 and 10,000 training steps. The training time took 282.19 seconds.
The exact solution and the PINN solution are shown in Fig. \ref{helmholtz1d-pinn}(a) and the pointwise error is shown in Fig. \ref{helmholtz1d-pinn}(b). The FFT of the exact solution and the PINN solution are shown in Fig. \ref{helmholtz1d-pinn}(c).
From the perspectives of pointwise errors and spectra, it can be seen that the accuracy of the PIMWNN method is higher than that of the PINN method. The  speed of the PIMWNN method is also faster than that of the PINN method.

\subsection{2D steady case}
\begin{figure}[!t]
\centering
	\subfigure[Exact solution.]{
		\includegraphics[width=0.3\textwidth]{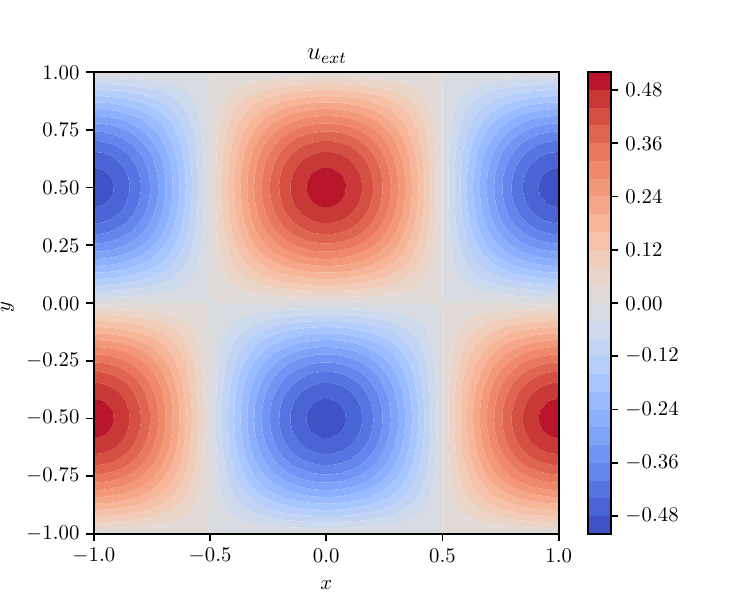}
	}
	\subfigure[PIMWNN solution.]{
		\includegraphics[width=0.3\textwidth]{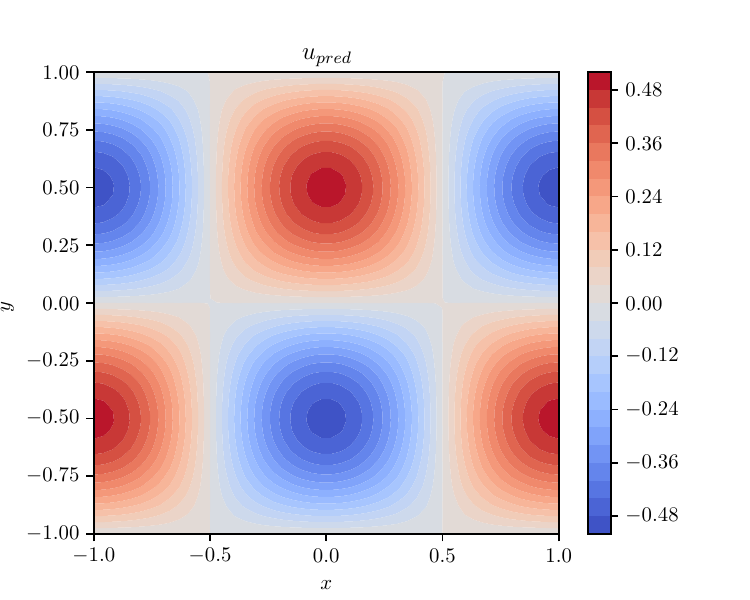}
	}	
	\subfigure[Pointwise error.]{
		\includegraphics[width=0.3\textwidth]{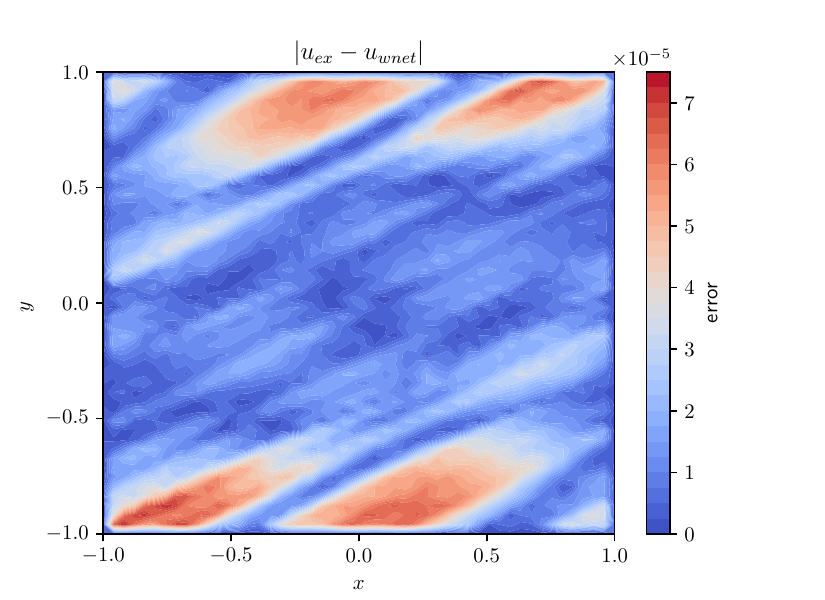}
	}	
\caption{Solution and error plots for 2D advection equation  Eq.(\ref{Tc-4}).}\label{tc-4}
\end{figure}

\begin{figure}[!t]
\centering
	\subfigure[Exact solution.]{
		\includegraphics[width=0.3\textwidth]{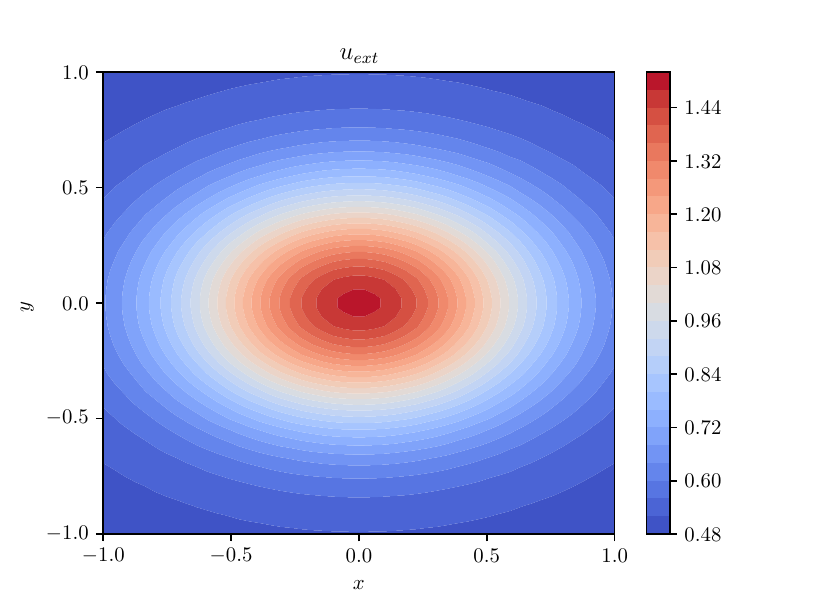}
	}
	\subfigure[PIMWNN solution.]{
		\includegraphics[width=0.3\textwidth]{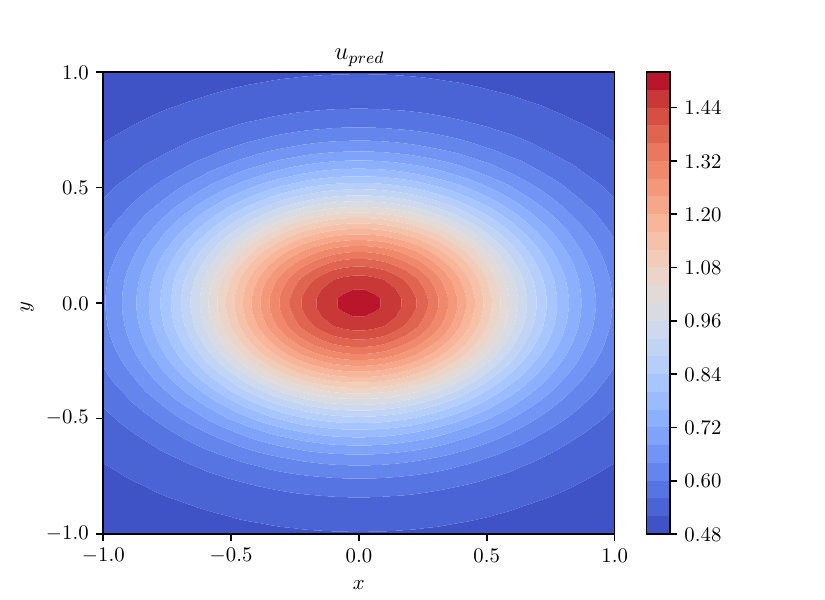}
	}	
	\subfigure[Pointwise error.]{
		\includegraphics[width=0.3\textwidth]{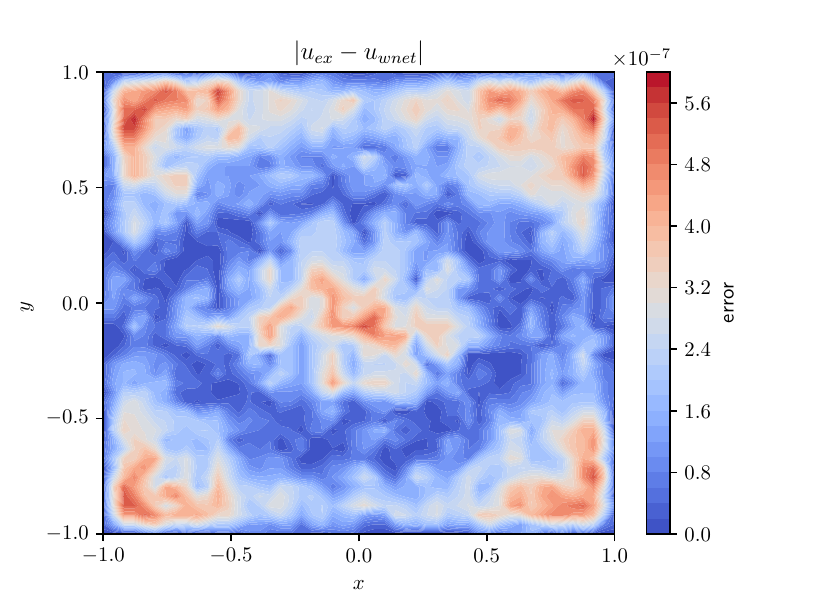}
	}	
\caption{Solution and error plots for 2D diffusion equation  Eq.(\ref{Tc-5}).}\label{tc-5}
\end{figure}

Next, we solve 2D stationary advection and diffusion equations, given by
\begin{subequations}
\begin{align}
au_{x} + bu_{y} = f, \quad (x,y) \in[-1,1]\times[-1,1],   \label{Tc-4} \\
(u_{xx} + u_{yy}) = f, \quad (x,y) \in[-1,1]\times[-1,1].   \label{Tc-5}
\end{align}
\end{subequations}
The Dirichlet boundary conditions and the function $f$ are chosen such that the exact solution $u$ are given by
\begin{subequations}
\begin{align}
u(x,y) &= \frac{1}{2}\cos(\pi x)\sin(\pi y), \\
u(x,y) &= \frac{1}{2}+\exp(-(2x^2+4y^2)), 
\end{align}
\end{subequations}
respectively. The results for these two cases are given in Fig. \ref{tc-4} and Fig. \ref{tc-5}, where the scale $J^x \times J^y = 4\times 4$ and $J^x \times J^y = 3\times 3$, respectively. The summary of the experiments of Eq. (\ref{Tc-4}) and Eq. (\ref{Tc-5}) are given in Table \ref{Tc4_5}.  We can observe again that the PIMWNN method can provide accurate numerical solutions of the 2D stationary advection and diffusion equations, and the relative $L^2$ errors of the PIMWNN solutions are reduced with the increase of the scale $J$. 

In \cite{Dwivedi2019PhysicsIE}, the unified deep ANN algorithm took 5500 points to achieve an order of accuracy of $10^{-5}$ using about 8 hours in these test cases. In comparison, we solve the case with 5000 points to achieve an order of accuracy of $10^{-4}$ and $10^{-7}$ taking only about 1 second, respectively.

\begin{table*}[!t]
\caption{Summary of experiments for test cases.}\label{Tc4_5}
\centering
\begin{tabular}{cccccccc}
\hline
  Equations & $N_f$ & $N_b$ & $J^x_0\times J^y_0$ & $J^x\times J^y$ & $N$ & $e_{L^2}$ & $t_{cpu}(s)$ \\
\hline
        & 5000 & 400 & 0$\times$0 & 2$\times$2 & 12$\times$12 & 1.235e-03 & 3.235e-02\\
Eq. (\ref{Tc-4}) & 5000 & 400 & 0$\times$0 & 3$\times$3 & 21$\times$21 & 2.227e-04 & 1.174e-01\\
        & 5000 & 400 & 0$\times$0 & 4$\times$4 & 38$\times$38 & 1.074e-04 & 1.077e-00\\
\hline
         & 1000 & 100 & 0$\times$0 & 1$\times$1 &  7$\times$7   & 4.984e-03  & 4.179e-01\\
Eq. (\ref{Tc-5}) & 1000 & 100 & 0$\times$0 & 2$\times$2 & 12$\times$12  & 5.117e-06  & 9.938e-01\\
       & 1000 & 100 & 0$\times$0 & 3$\times$3 & 21$\times$21  & 3.056e-07  & 3.943e-00\\
\hline
Eq. (\ref{Tc-9})    & 10000 & 100 & -8$\times$-4 & 14$\times$8 & 152$\times$46 & 7.322e-04 & 2.525e02\\
Eq. (\ref{Tc-10})  & 5000 & 100 &  -7$\times$-2 & 11$\times$4 &  52$\times$16 & 3.413e-04 & 4.831e-01\\
\hline
\end{tabular}
\end{table*}

\subsection{2D diffusion equation on a complex domain}
In this example, we solve 2D diffusion equation on a complex domain, given by
\begin{equation}
	\begin{split}
         (u_{xx} +u_{yy}) &= f, \quad (x,y) \in \Omega  \label{Tc-6}\\
	\end{split}
\end{equation}
where $\Omega = \{(x,y)| x=0.5 + r(\theta)\cos(\theta), y=0.5 + r(\theta)\sin(\theta), r(\theta) = 0.2+0.15\sin(5\theta), \theta \in [0,2\pi) \}$, which has a complex boundary,  something like a flower. The function $f$ is chosen such that the exact solution is given by
\begin{equation}
	\begin{split}
		u(x,y) &= 16x(1-x)y(1-y)\\
		&(0.5 + \frac{\arctan (20(0.25^2-(x-0.5)^2-(y-0.5)^2) )}{\pi}). 
	\end{split}
\end{equation}
The Dirichlet boundary condition is given by
\begin{equation}
	\begin{split}
		&u(x,y)|_\Gamma = 16x(1-x)y(1-y)\cdot \\
		&\left(0.5 + \frac{\arctan (20(0.25^2-(x-0.5)^2-(y-0.5)^2) )}{\pi}\right)|_\Gamma, \\
	\end{split}
\end{equation}
where $\Gamma$ denotes the complex boundary of domain $\Omega$. 

\begin{figure}[!t]
\centering
	\subfigure[Exact solution.]{
		\includegraphics[width=0.3\textwidth]{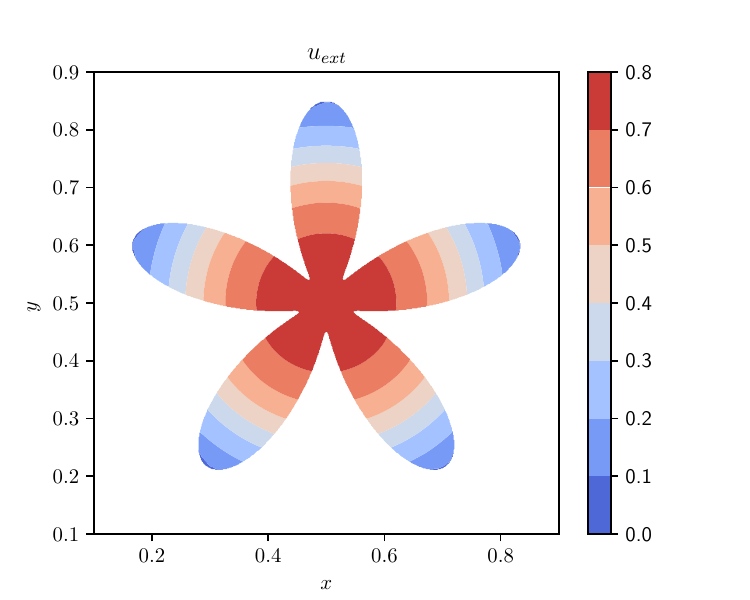}
	}
	\subfigure[PIMWNN solution.]{
		\includegraphics[width=0.3\textwidth]{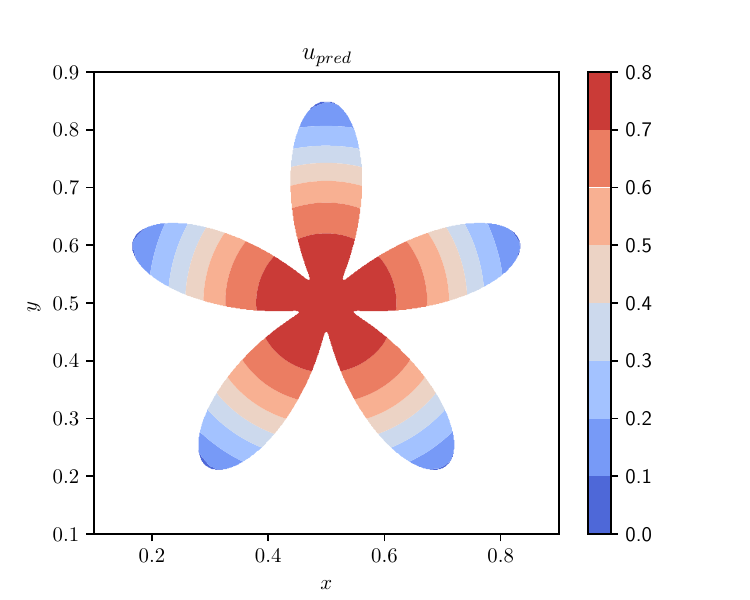}
	}	
	\subfigure[Pointwise error.]{
		\includegraphics[width=0.3\textwidth]{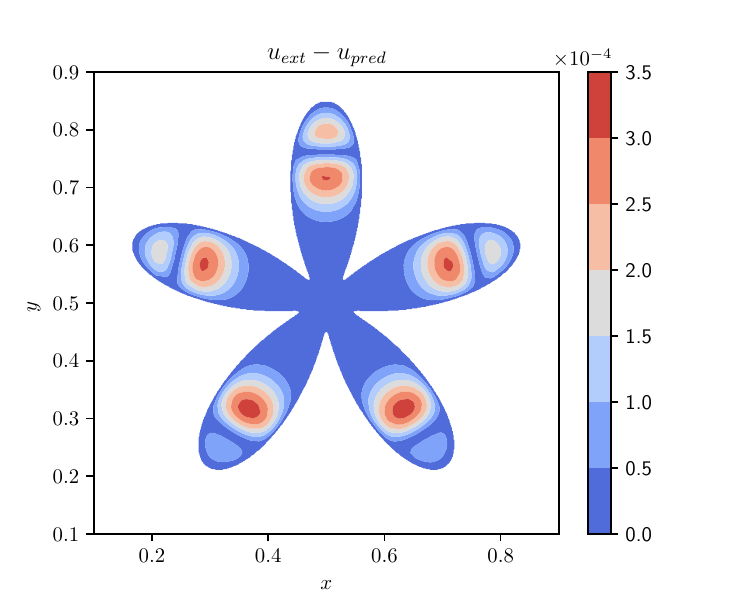}
	}	
\caption{Solution and error plots for 2D diffusion equation  Eq.(\ref{Tc-6}) on a complex domain.}\label{tc-6-1}
\end{figure}

The exact solution is shown in Fig. \ref{tc-6-1}(a) and the PIMWNN solution  is Fig. \ref{tc-6-1}(b), where the scale $J^x \times J^y = 4\times 4$. The pointwise error is shown in Fig. \ref{tc-6-1}(c). We can observe that the PIMWNN method can provide accurate numerical solutions of the 2D diffusion equation on a complex domain.

\subsection{Unsteady cases}

\subsubsection{Linear advection equation with complicated initial condition}
\begin{figure}[!t]
\centering
	\subfigure[Exact and PINN solution of linear advection of a high frequency wave packet. Red: PINN, Blue: exact.\cite{Dwivedi2019PhysicsIE}]{
        \includegraphics[width=0.3\textwidth]{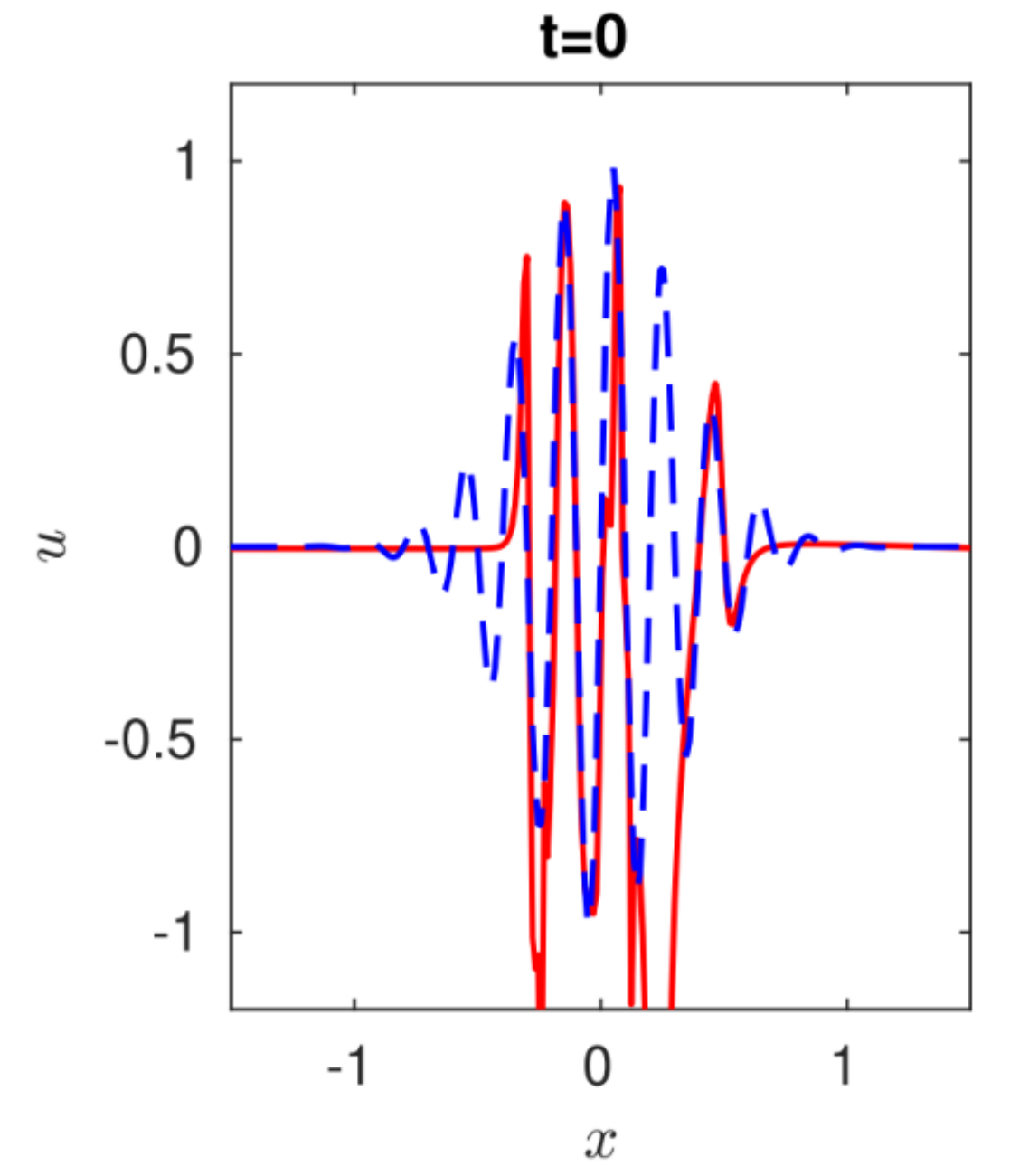}
        \includegraphics[width=0.3\textwidth]{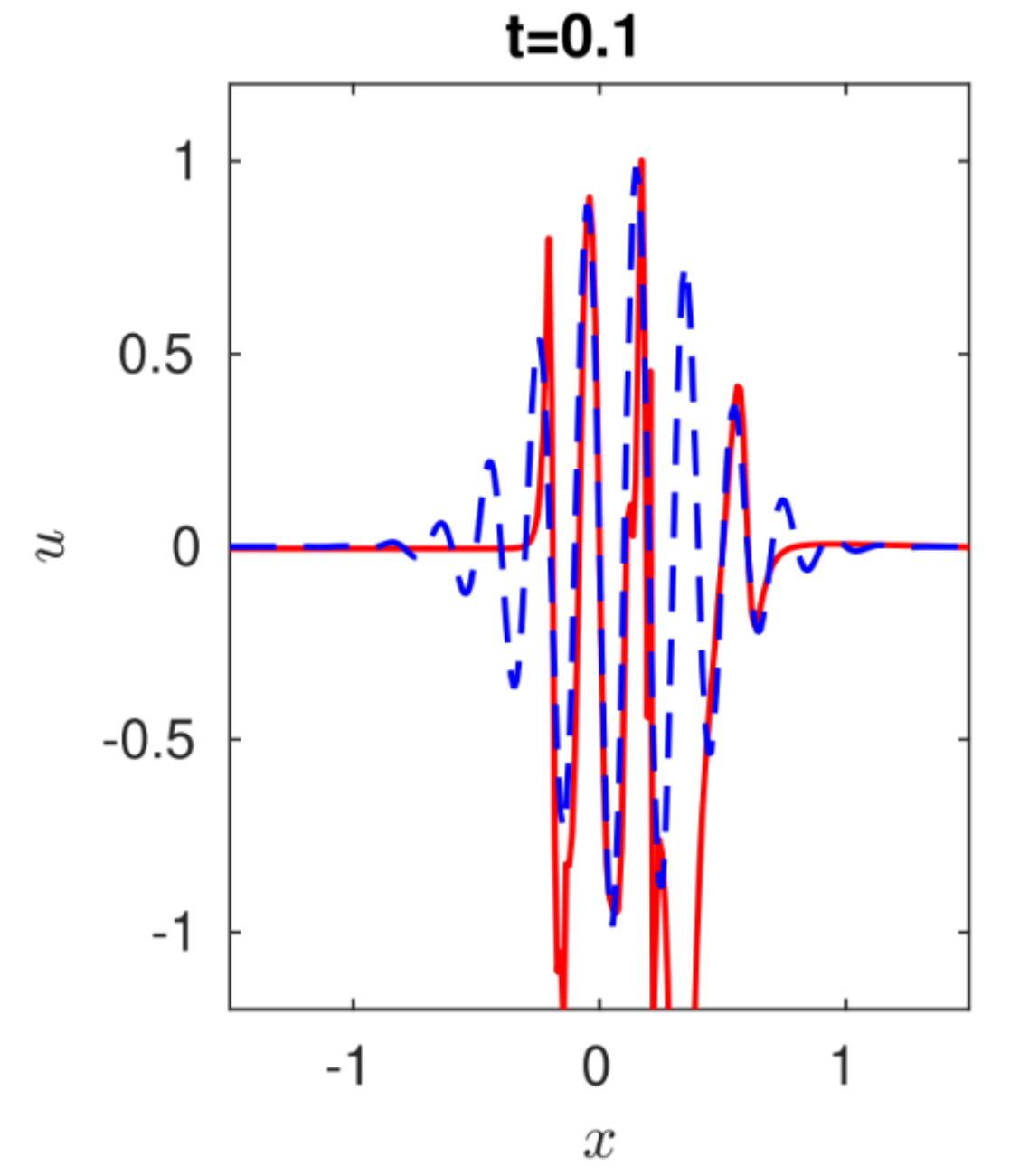}
        \includegraphics[width=0.3\textwidth]{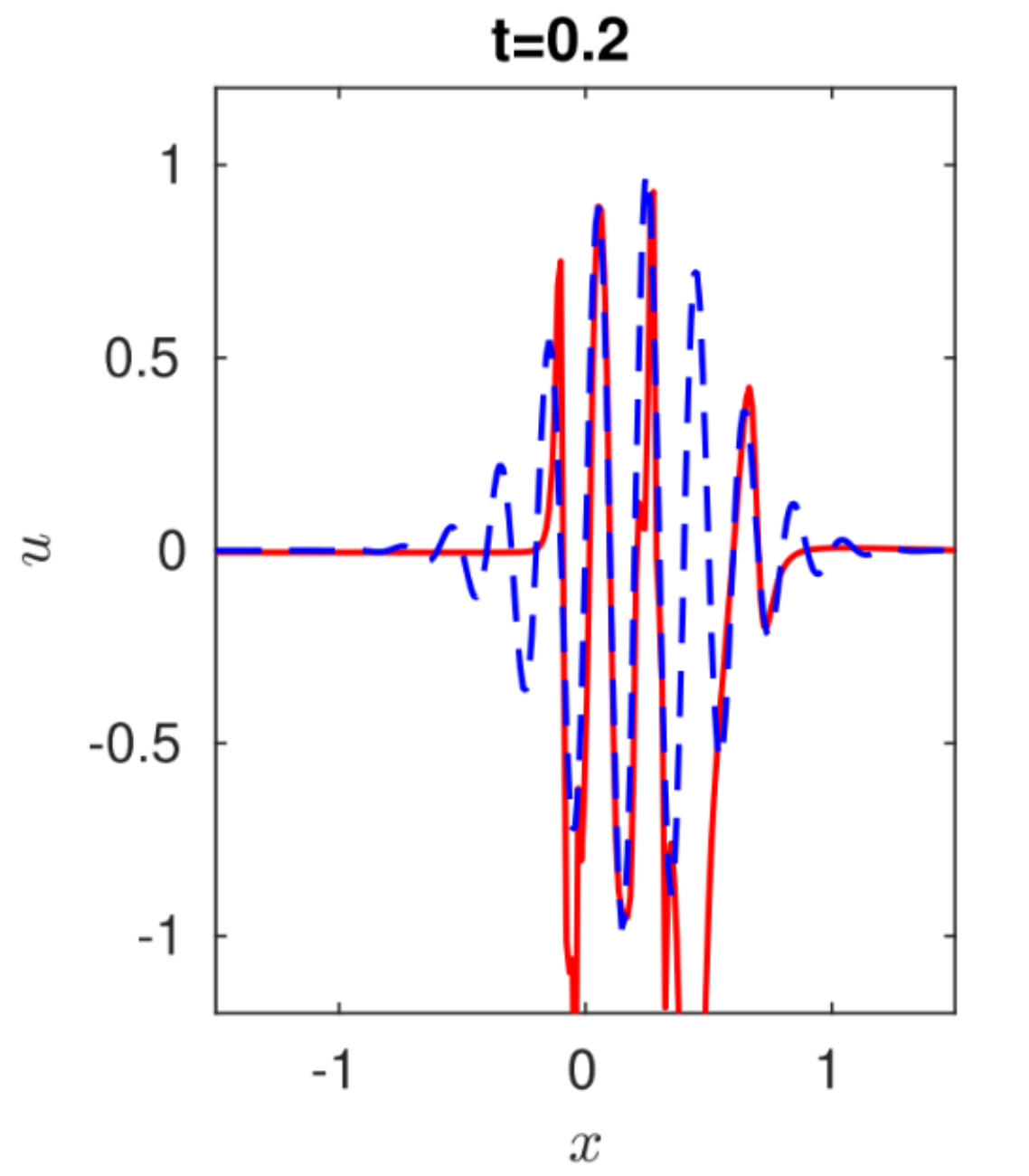}
	}
	\subfigure[Exact and PIMWNN solutions of linear advection of a high frequency wave packet. Red: PIMWNN, Blue: exact.]{
        \includegraphics[width=2.in]{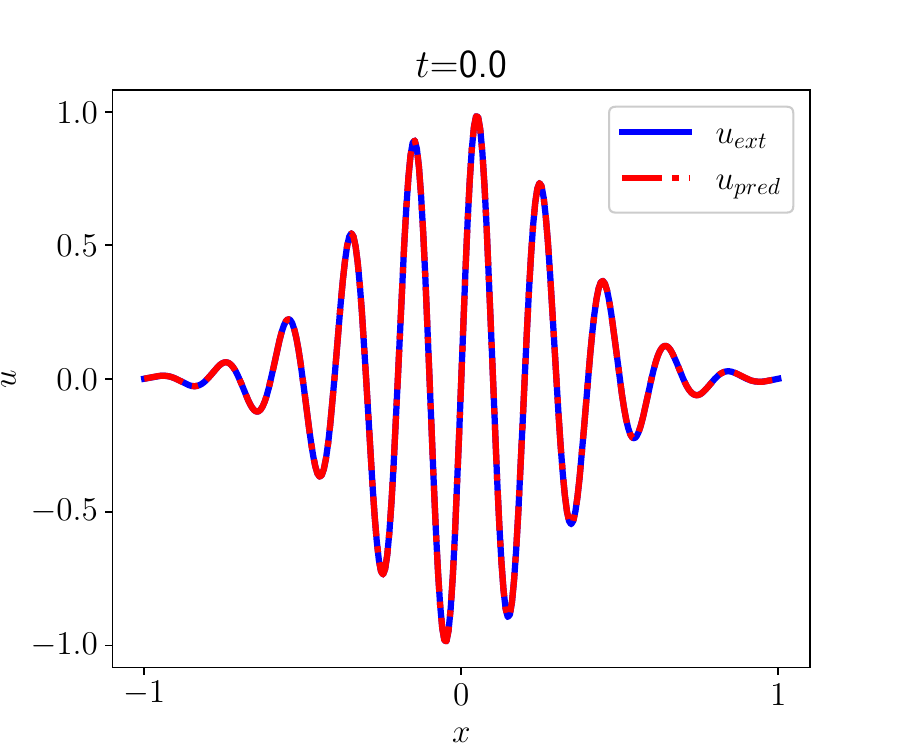}
        \includegraphics[width=2.in]{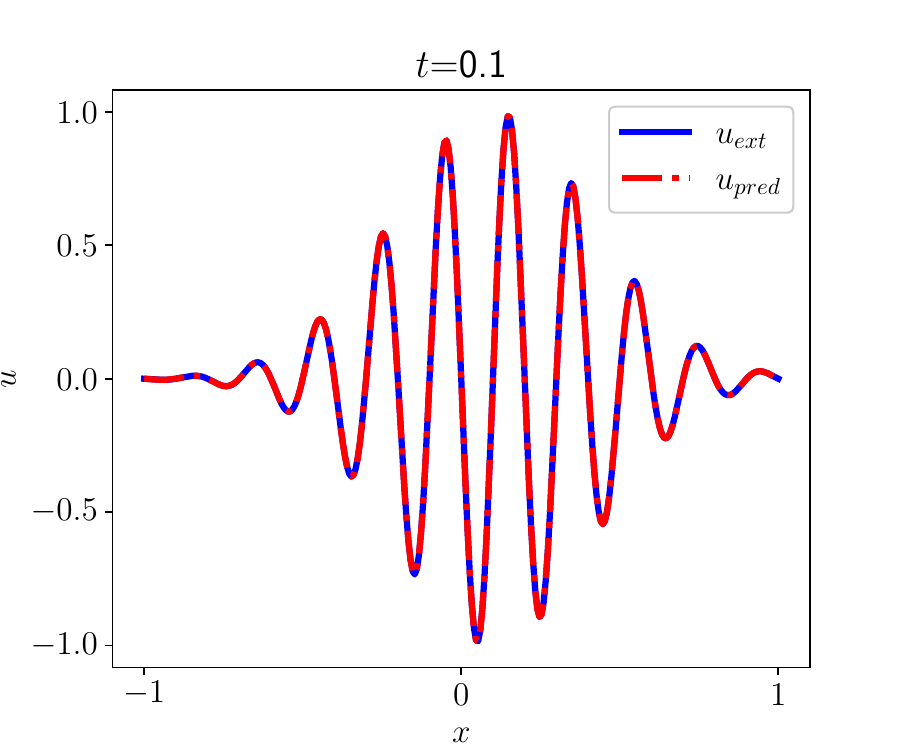}
        \includegraphics[width=2.in]{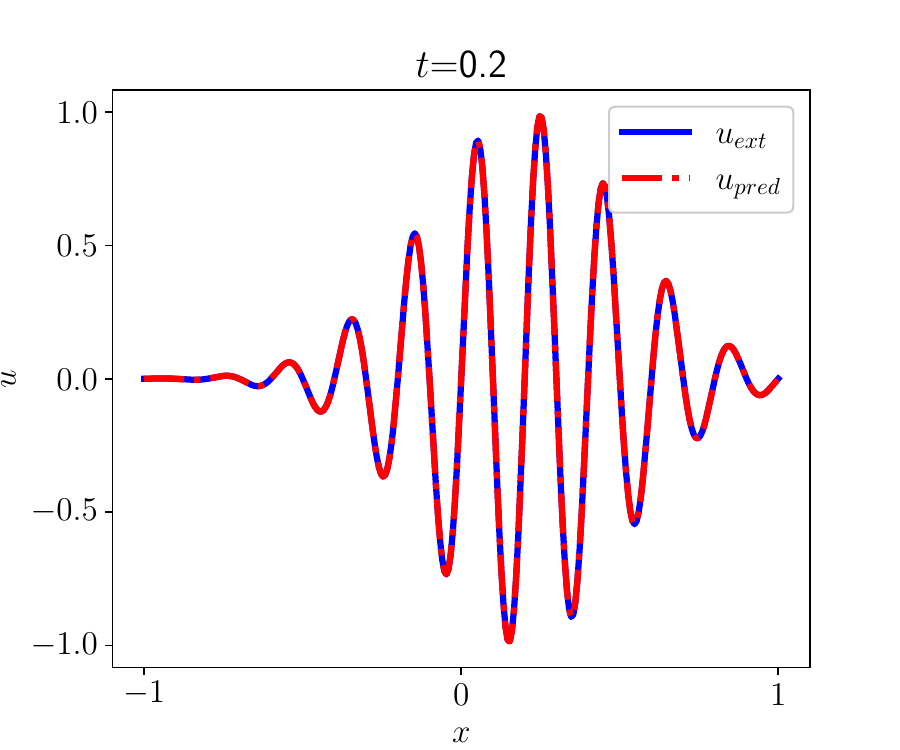}
	}	
	\subfigure[Pointwise errors of PIMWNN solution at different times.]{
        \includegraphics[width=2.in]{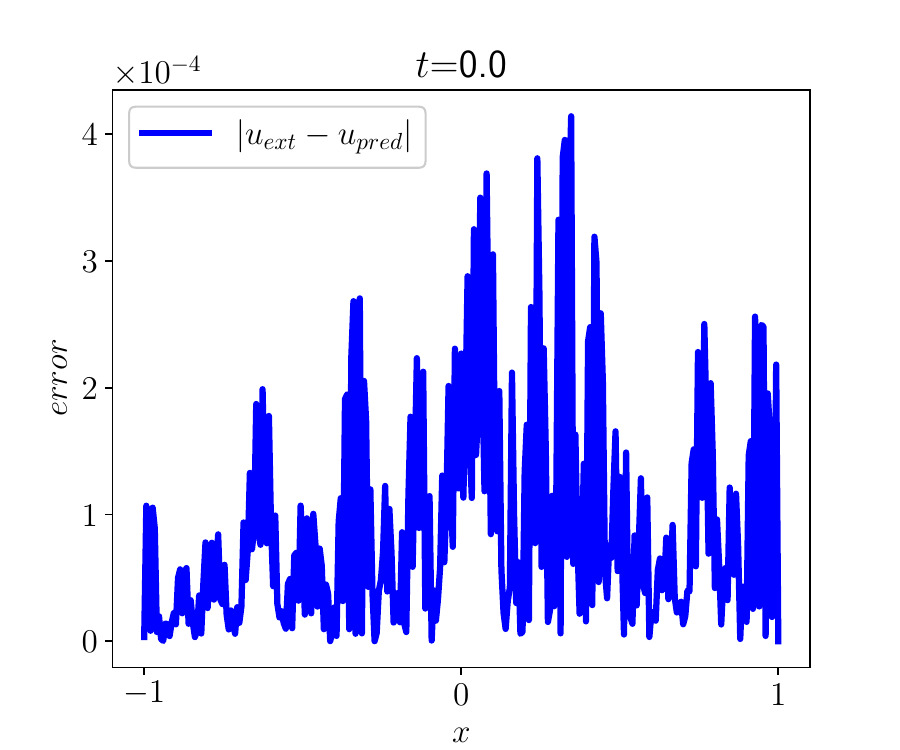}
        \includegraphics[width=2.in]{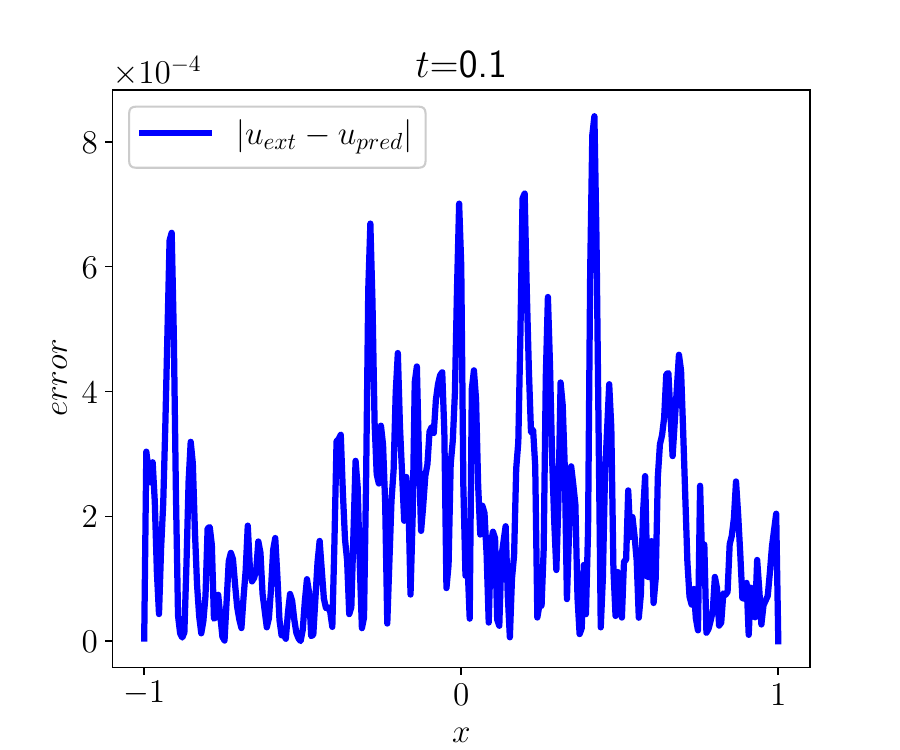}       
        \includegraphics[width=2.in]{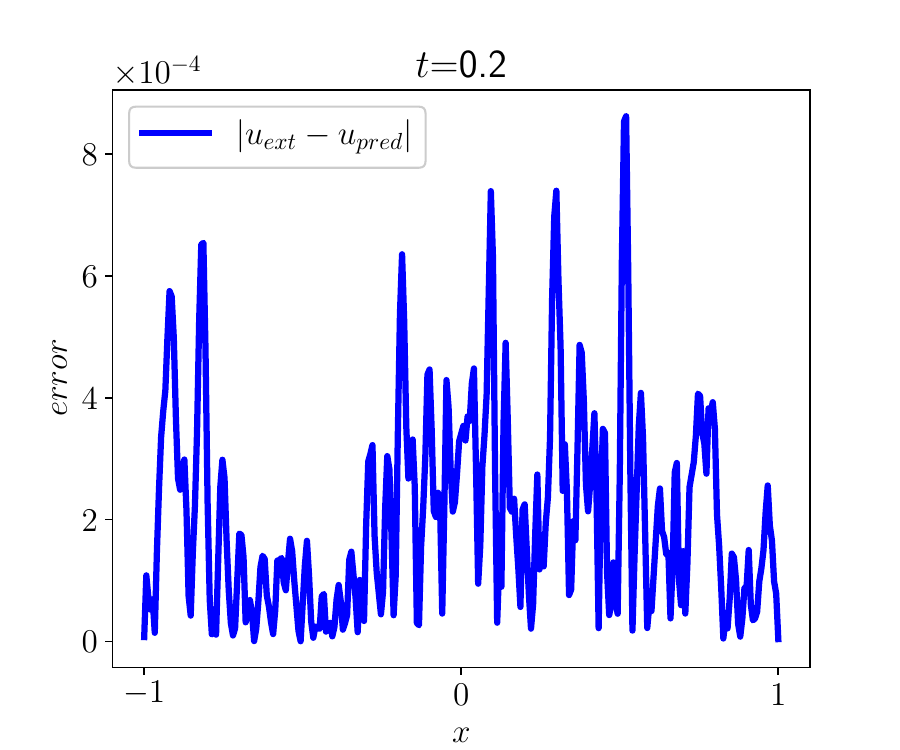}
	}	
\caption{Solution and error plots for 1D unsteady advection Eq.(\ref{Tc-9}).}\label{tc-8}
\end{figure}

In this example, we solve a 1D unsteady linear advection equation with a complicated initial condition, given by
\begin{equation}
\begin{split}
&u_t + u_x = 0, \quad (x,t)\in (-1,1)\times (0,0.5],  \label{Tc-9} \\
&u(x,0) = \exp(-5x^2)\sin(10\pi x), \quad x\in [-1,1], \\
&u(-1,t) = u(1,t) = 0, \quad t \in (0,0.5].
\end{split}
\end{equation}
The exact solution is given by
\begin{equation}
\begin{split}
u(x,t) &= \exp(-5(x-t)^2)\sin(10\pi (x-t)).
\end{split}
\end{equation}
This example is taken from \cite{Dwivedi2019PhysicsIE}. The author solved this equation with a deep PINN. The deep PINN consists of 9 hidden layers with 20 neurons each. The exact and PINN solutions are shown in Fig. \ref{tc-8}(a), while the PIMWNN solutions and the corresponding pointwise errors are shown in Figs. \ref{tc-8}(b) and \ref{tc-8}(c). 

\subsubsection{Linear advection equation with a sharp-peaked Gaussian initial condition}
\begin{figure}[!t]
\centering
	\subfigure[Exact and PIELM solutions of linear advection of a sharp-peaked Gaussian initial condition.  \cite{Dwivedi2019PhysicsIE}.]{
		\includegraphics[width=0.8\textwidth]{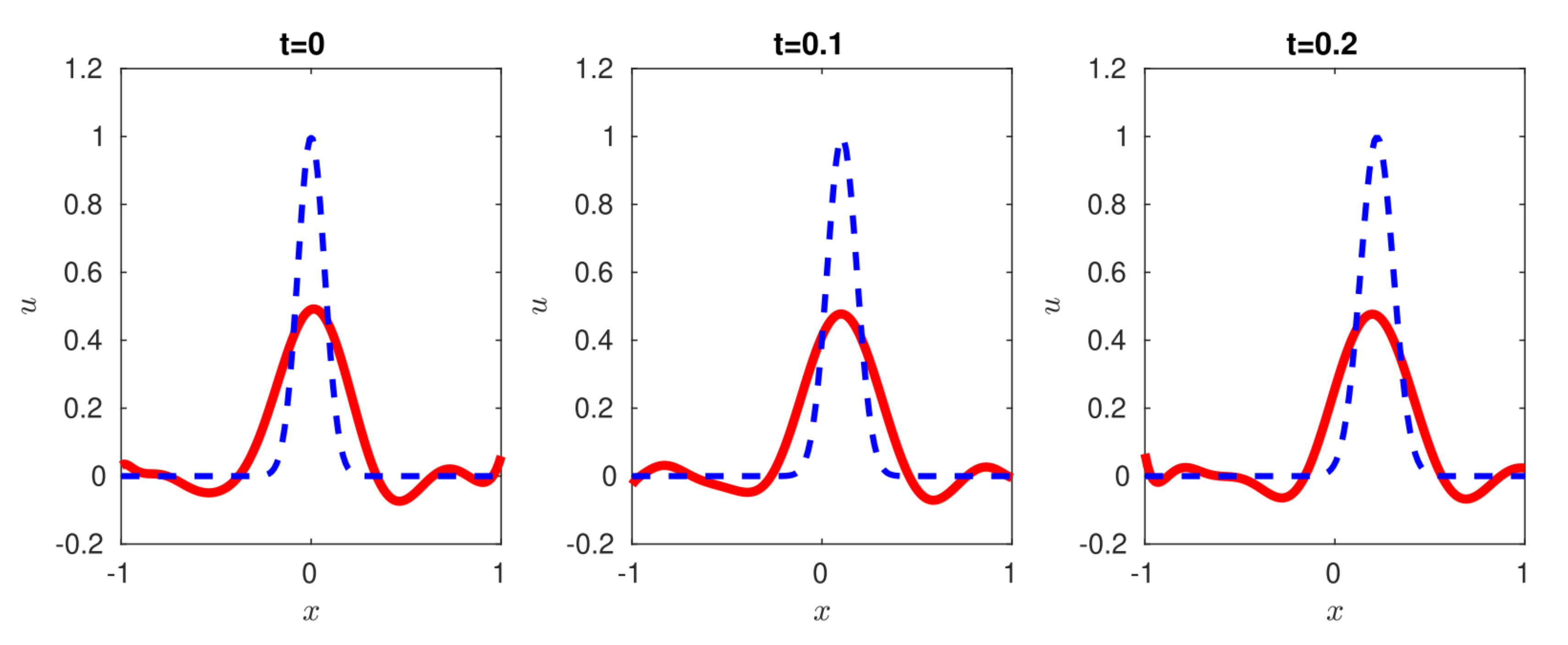}
	}
	
	\subfigure[Exact and PIMWNN solutions of linear advection of a sharp-peaked Gaussian initial condition. Red: PIMWNN, Blue: exact.]{
		\includegraphics[width=0.9\textwidth]{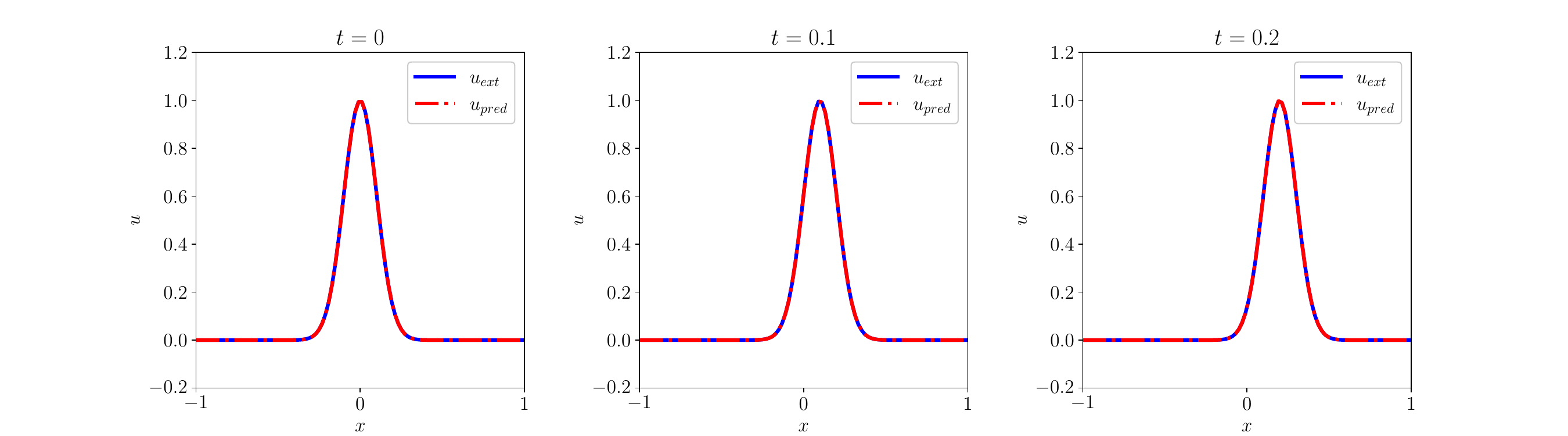}
	}	
	
	\subfigure[Pointwise errors of PIMWNN solution at different times.]{
		\includegraphics[width=0.9\textwidth]{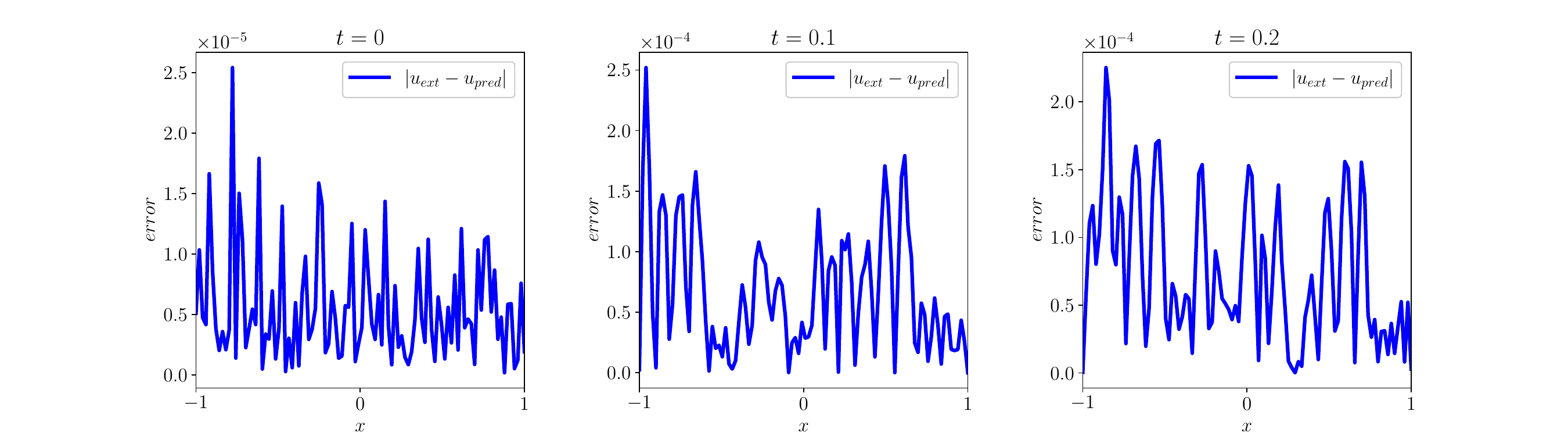}
	}	
\caption{Solution and error plots for 1D unsteady advection Eq.(\ref{Tc-10}).}\label{tc-9}
\end{figure}

In this example, we solve a 1D unsteady linear advection equation with a sharp-peaked Gaussian initial condition,
\begin{equation}
\begin{split}
&u_t + u_x = 0, \quad (x,t)\in (-1,1)\times (0,0.5], \label{Tc-10} \\
&u(x,0) = \exp(-50x^2), \quad x\in [-1,1], \\
&u(-1,t) = u(1,t) = 0, \quad t \in (0,0.5].
\end{split}
\end{equation}
The exact solution is given by
\begin{equation}
\begin{split}
u(x,t) &= \exp(-50(x-t)^2).
\end{split}
\end{equation}
This example is taken from \cite{Dwivedi2019PhysicsIE}. The authors solved this equation with the PIELM method. The exact and PIELM solutions are shown in Fig. \ref{tc-9}(a) , highlighting the PIELM algorithm's inability to effectively solve the problem. The PIMWNN solutions and the associated pointwise error are shown in Figs. \ref{tc-9}(b) and \ref{tc-9}(c), demonstrating the method's capability to solve the equation with high precision.

\subsubsection{Linear diffusion equation with the amplitude of solution grows exponentially over time}
To evaluate the performance of PIMWNN in a situation where the solution's amplitude grows exponentially over time, we solve the following example\cite{HUANG2025113676}:
\begin{equation}
\begin{split}
u_t - u_{xx} &= f(x,t) \ , \quad (x,t)\in (-1,1)\times (0,1], \label{unsteady-diffusion-eq} \\
u(x,0) &= \sin(k_1 x) - \cos(k_2x), \quad x \in [-1,1],\\
u(-1,t) &= u(1,t) = -e^t, \quad  t\in [0,1].
\end{split}
\end{equation}
The source term is given by 
\begin{equation}
f(x,t) = e^t((1+k_1^2)\sin(k_1 x) + (1+k_2^2)\cos(k_2x)),
\end{equation}

\begin{figure}[!t]
\centering
	\subfigure[Exact solution.]{
		\includegraphics[width=0.3\textwidth]{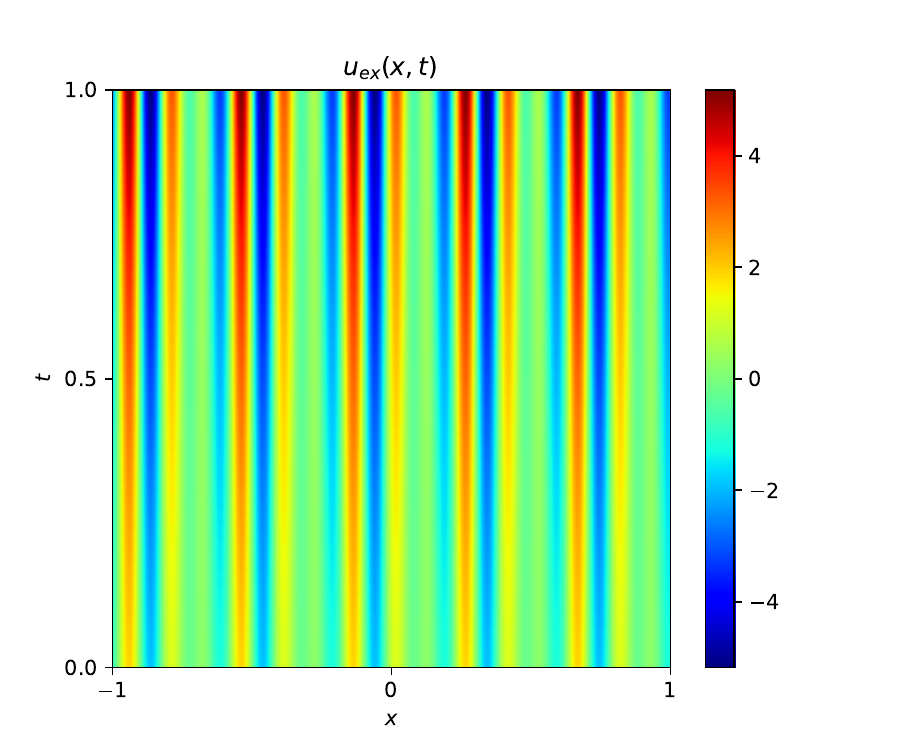}
	}
	\subfigure[PIMWNN solution.]{
		\includegraphics[width=0.3\textwidth]{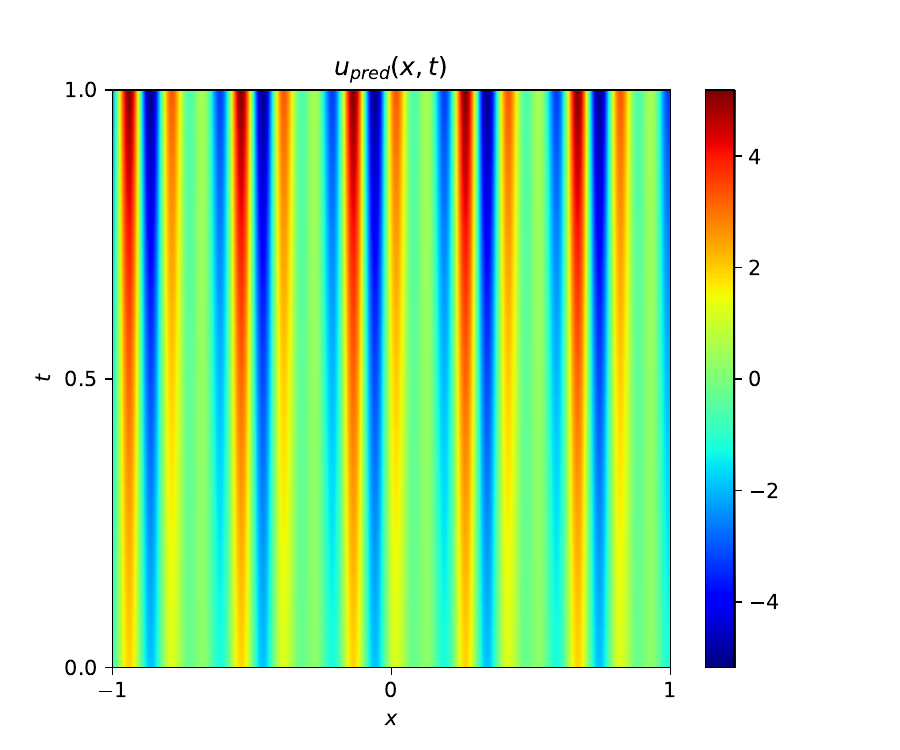}
	}	
	\subfigure[PIMWNN pointwise error.]{
		\includegraphics[width=0.3\textwidth]{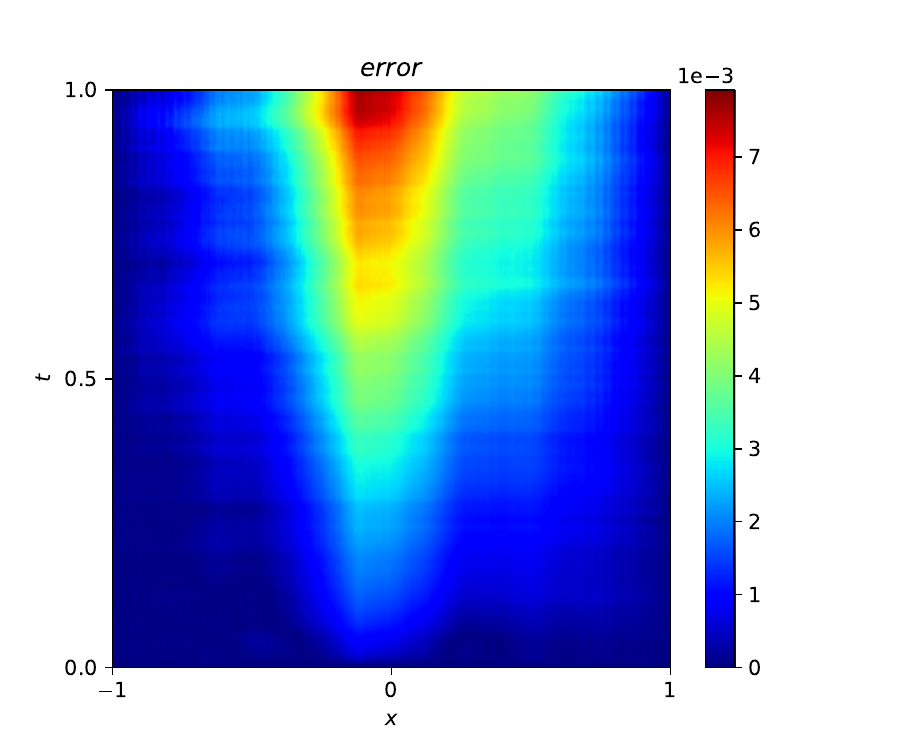}
	}	
    \subfigure[PINN solution\cite{HUANG2025113676}.]{
		\includegraphics[width=0.4\textwidth]{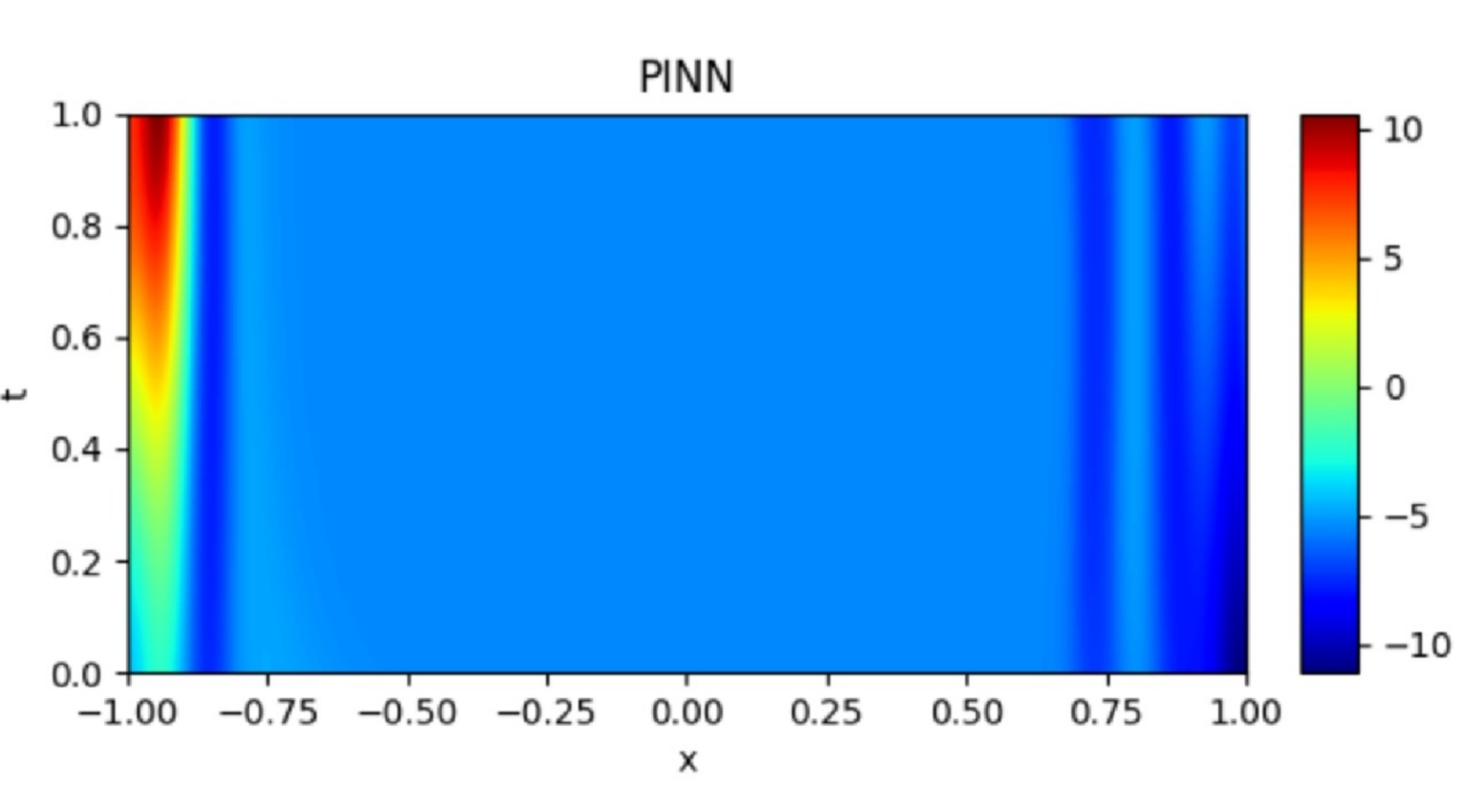}
	}	
	\subfigure[PINN pointwise error\cite{HUANG2025113676}.]{
		\includegraphics[width=0.4\textwidth]{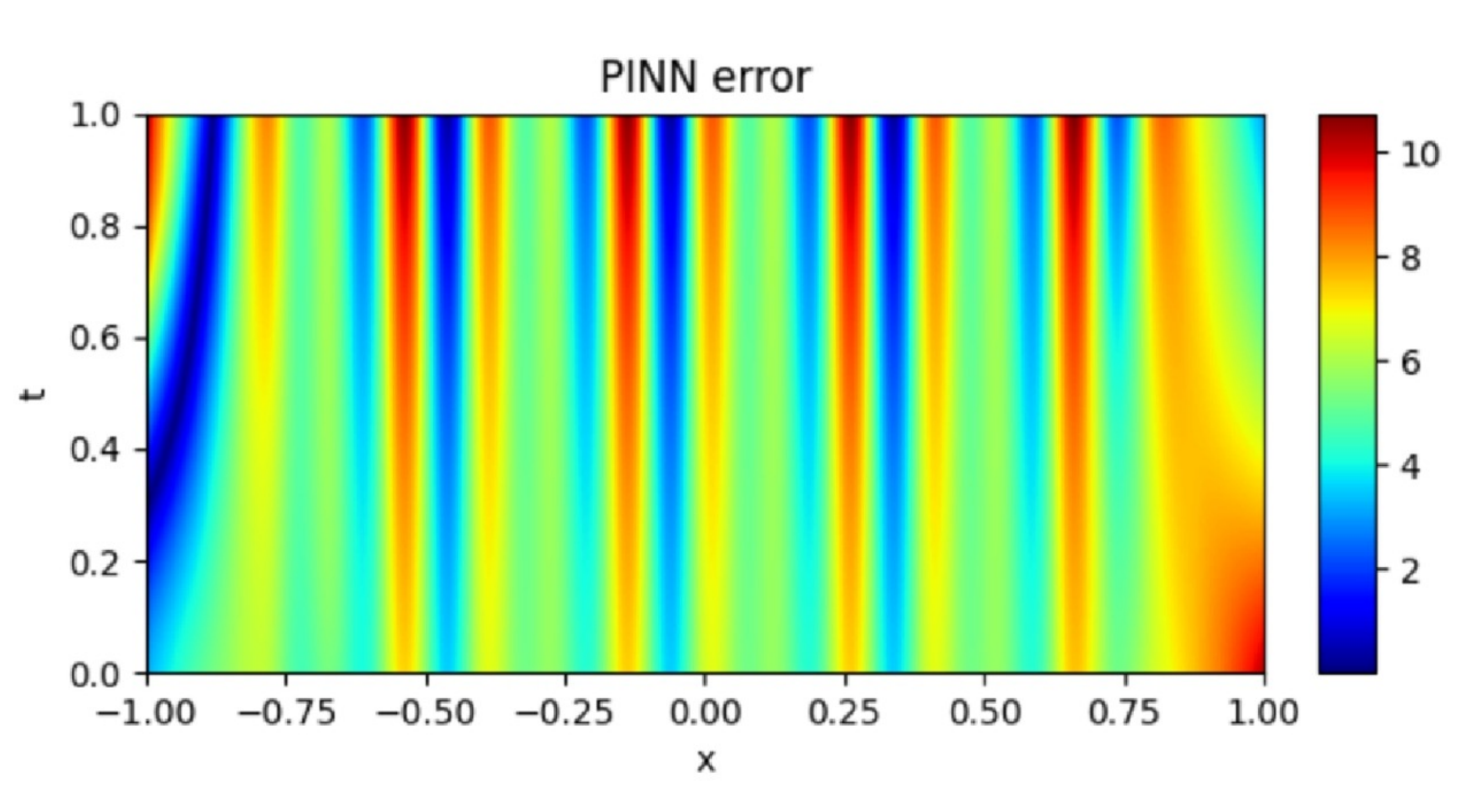}
	}	
\caption{Solution and error plots for diffusion equation  Eq.(\ref{unsteady-diffusion-eq}).}\label{unsteady-diffusion-fig}
\end{figure}

The exact solution and the PIMWNN solution of the above problem are shown in Fig. \ref{unsteady-diffusion-fig}(a) and Fig. \ref{unsteady-diffusion-fig}(b). the pointwise error is shown in Fig. \ref{unsteady-diffusion-fig}(c). The PINN solution and pointwise error shown in \ref{unsteady-diffusion-fig}(d) and \ref{unsteady-diffusion-fig}(e) are given in reference\cite{HUANG2025113676}.
We can see that PIMWNN solution matches well with the exact solution. The pointwise error shown in Fig. \ref{unsteady-diffusion-fig}(c) is of the order of $10^{-3}$. However, the PINN method fails to solve the problem correctly.

\subsection{Burger's equation}
The Burgers' equation arises in various areas of applied mathematics, including fluid mechanics, nonlinear acoustics, and gas dynamics, etc. This equation is a nonlinear partial differential equation that simulates the propagation and reflection of shock waves. In this example, we solve the Burger's equation along with Dirichlet boundary conditions,
\begin{equation}
\begin{split}
u_t + uu_x &= \frac{0.01}{\pi} u_{xx} \ , \quad (x,t)\in (-1,1)\times (0,1], \label{Tc-12} \\
u(x,0) &= -\sin(\pi x), \\
u(-1,t) &= u(1,t) = 0.
\end{split}
\end{equation}
This example is taken from \cite{raissi2019}. The exact solution is given in Package DeepXDE \cite{Lu2021deepxde}.

\begin{figure}[!t]
\centering
	\subfigure[Exact solution.]{
		\includegraphics[width=0.3\textwidth]{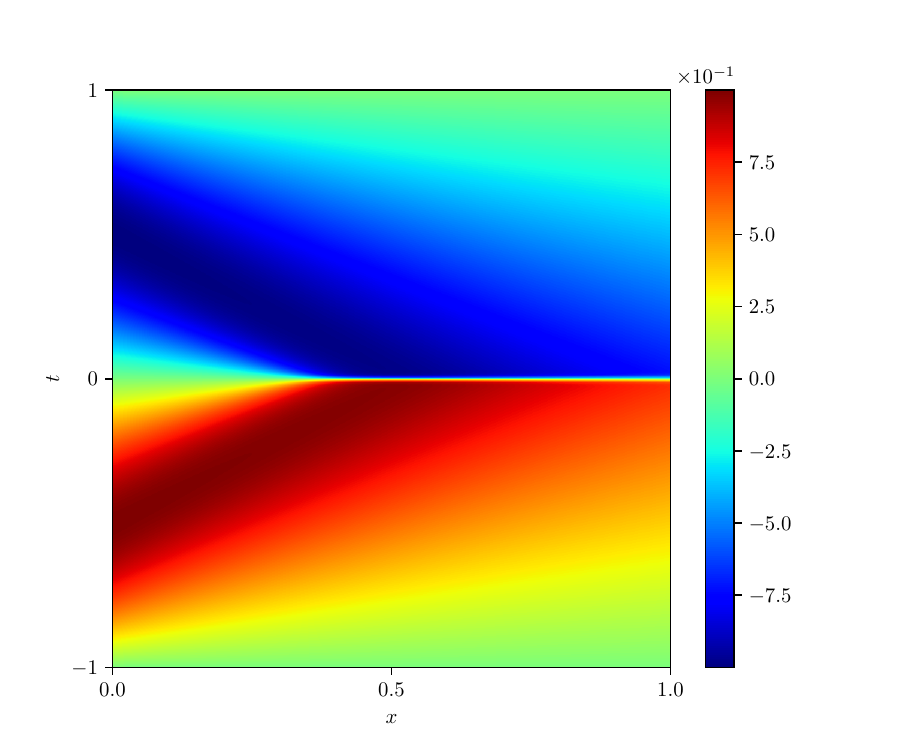}
	}
	\subfigure[PIMWNN solution.]{
		\includegraphics[width=0.3\textwidth]{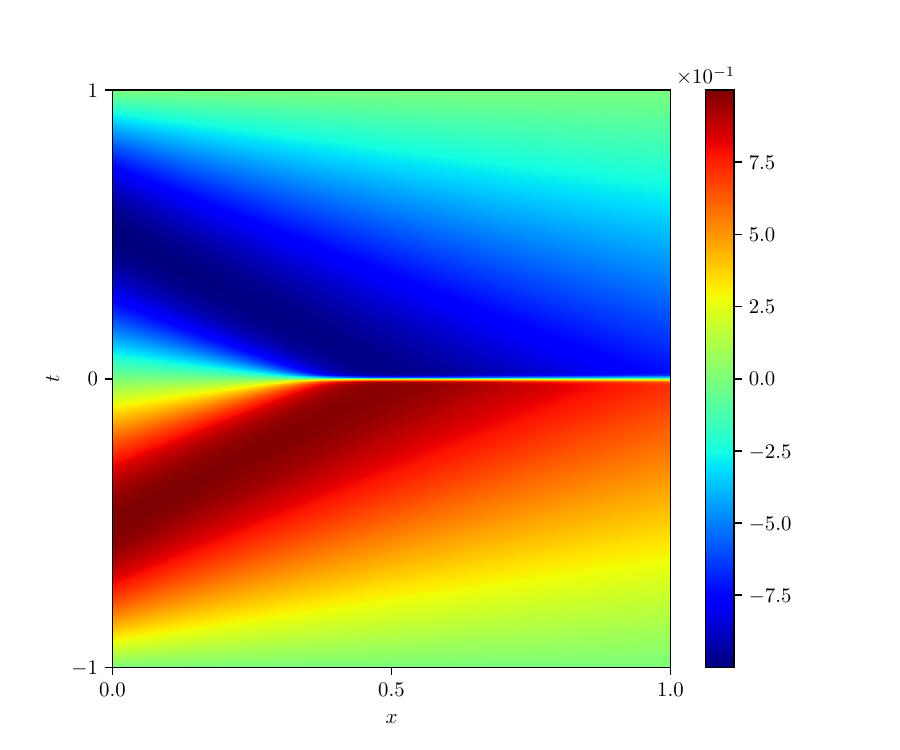}
	}	
	\subfigure[PIMWNN pointwise error.]{
		\includegraphics[width=0.3\textwidth]{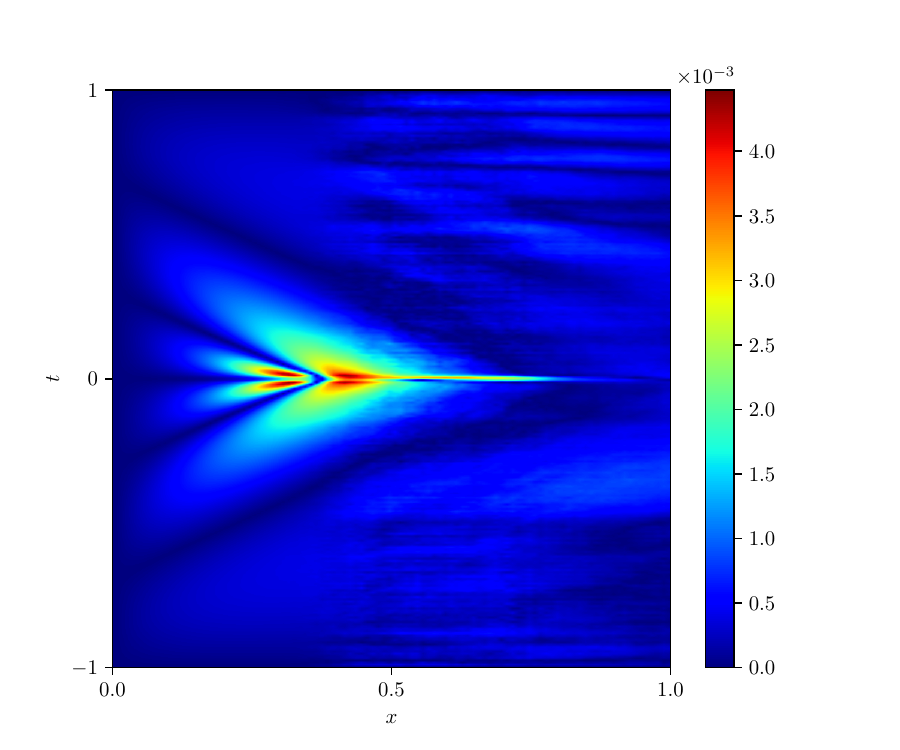}
	}	
	\subfigure[RNN solution.]{
		\includegraphics[width=0.3\textwidth]{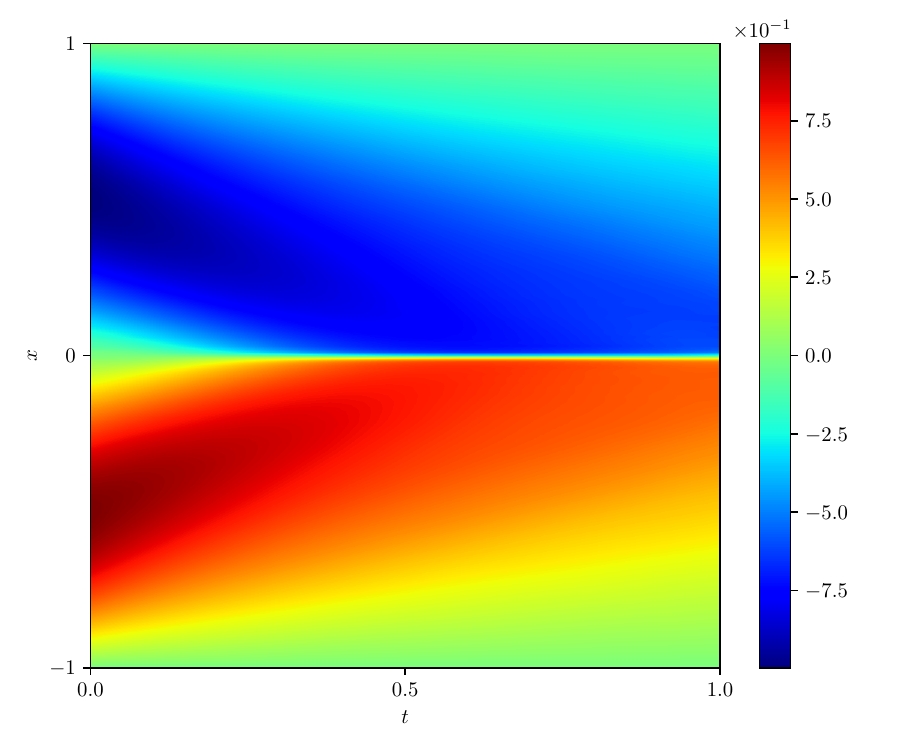}
	}	
	\subfigure[RNN pointwise error.]{
		\includegraphics[width=0.3\textwidth]{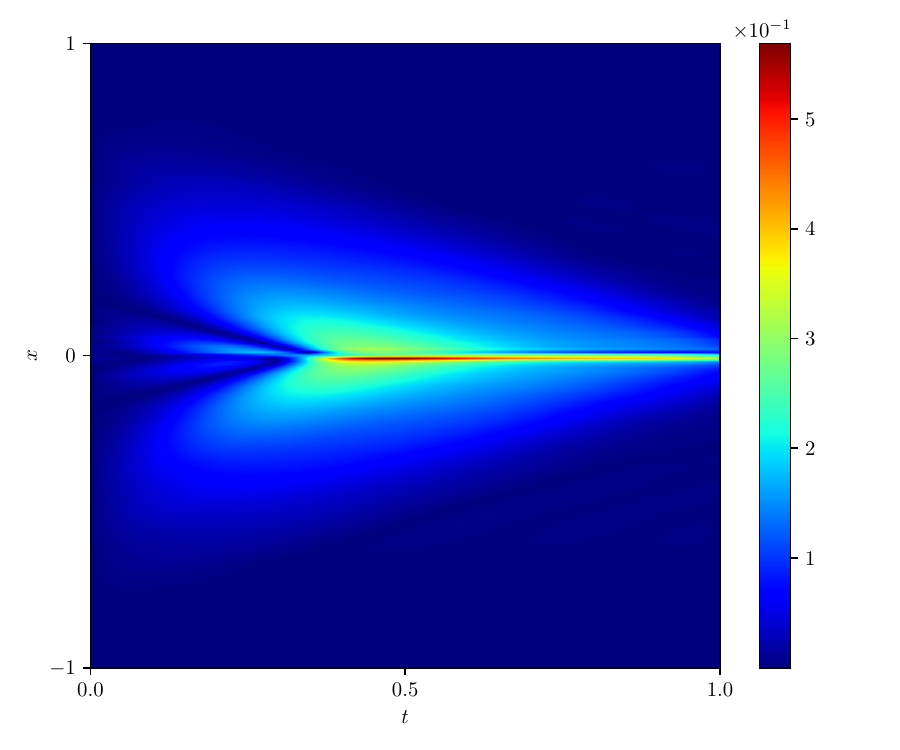}
	}	
\caption{Solution and error plots for Burger's equation  Eq.(\ref{Tc-12}).}\label{burgers-fig}
\end{figure}

\begin{figure}[!t]
\centering
	\subfigure[Exact and PIMWNN solutions of Burger's equation.]{
        \includegraphics[width=0.3\textwidth]{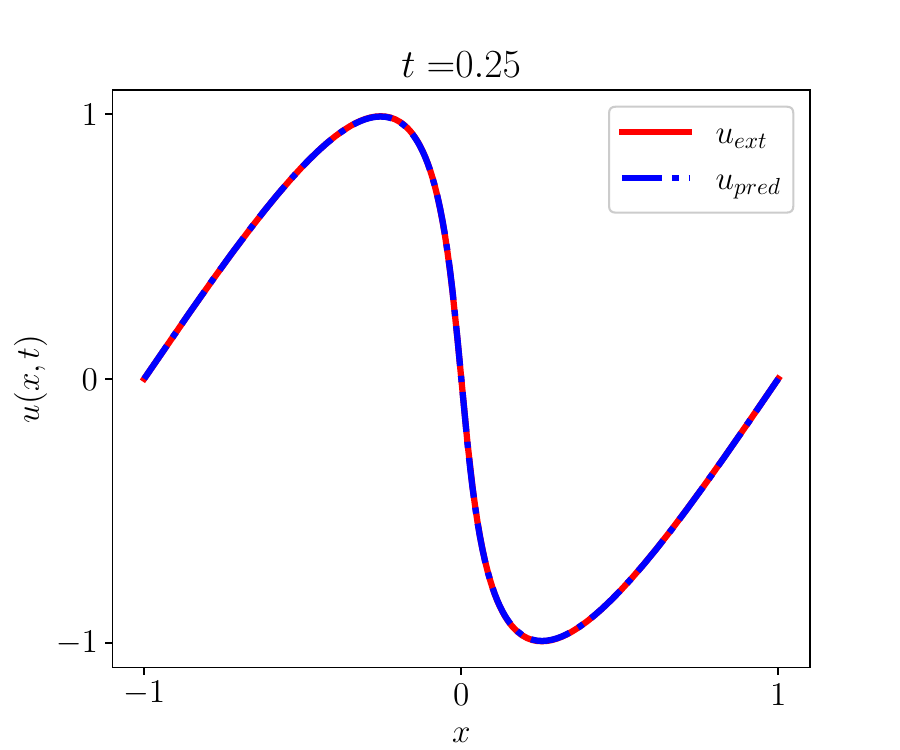}
        \includegraphics[width=0.3\textwidth]{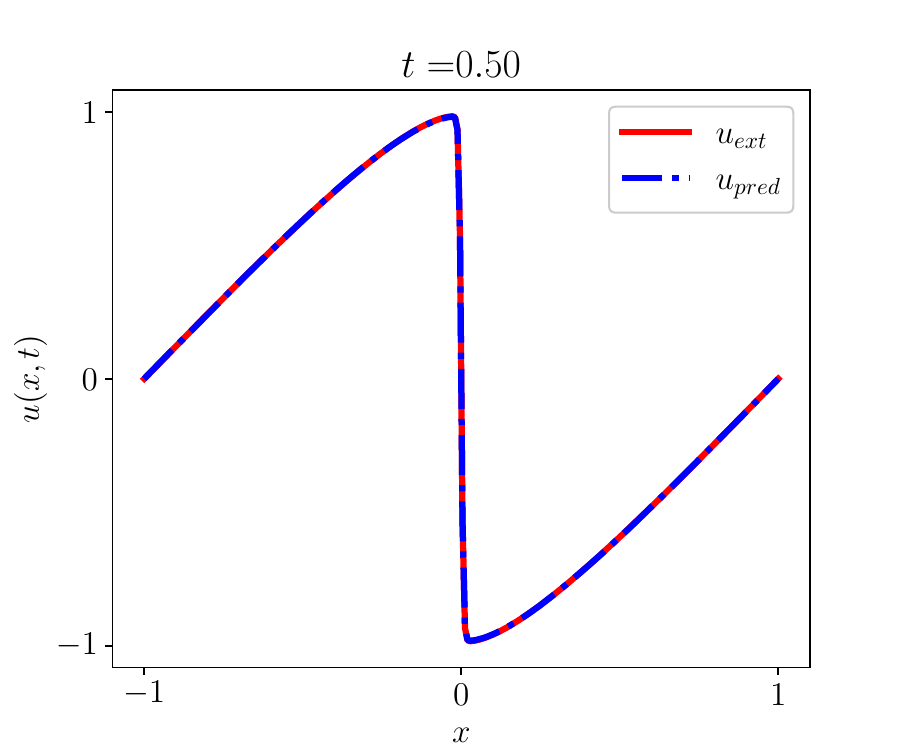}
        \includegraphics[width=0.3\textwidth]{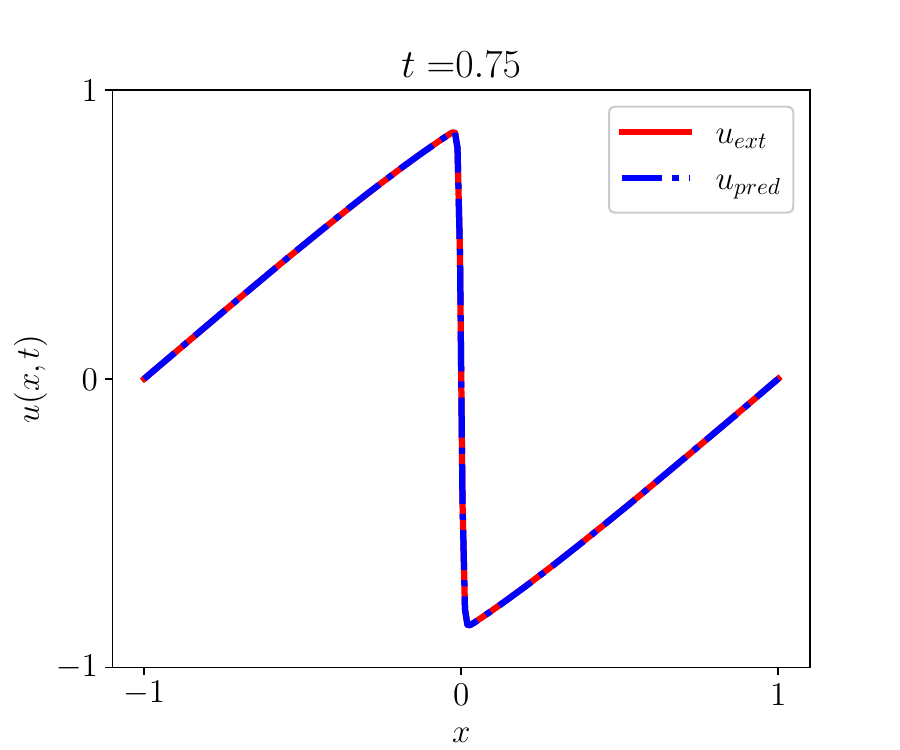}
	}
	\subfigure[Pointwise error of PIMWNN solutions of Burger's equation.]{
        \includegraphics[width=2.in]{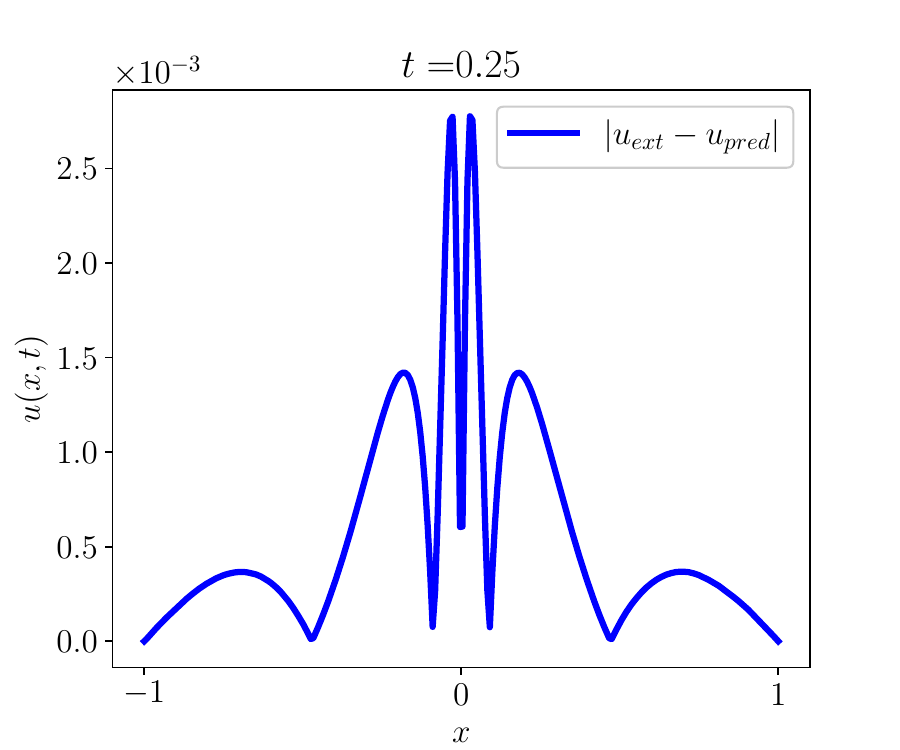}
        \includegraphics[width=2.in]{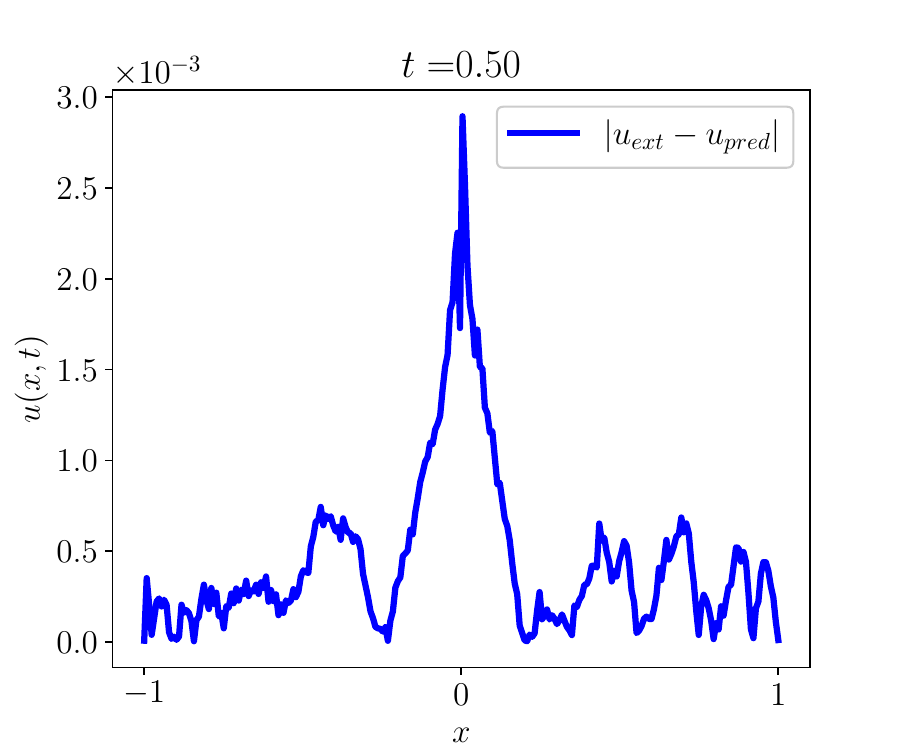}
        \includegraphics[width=2.in]{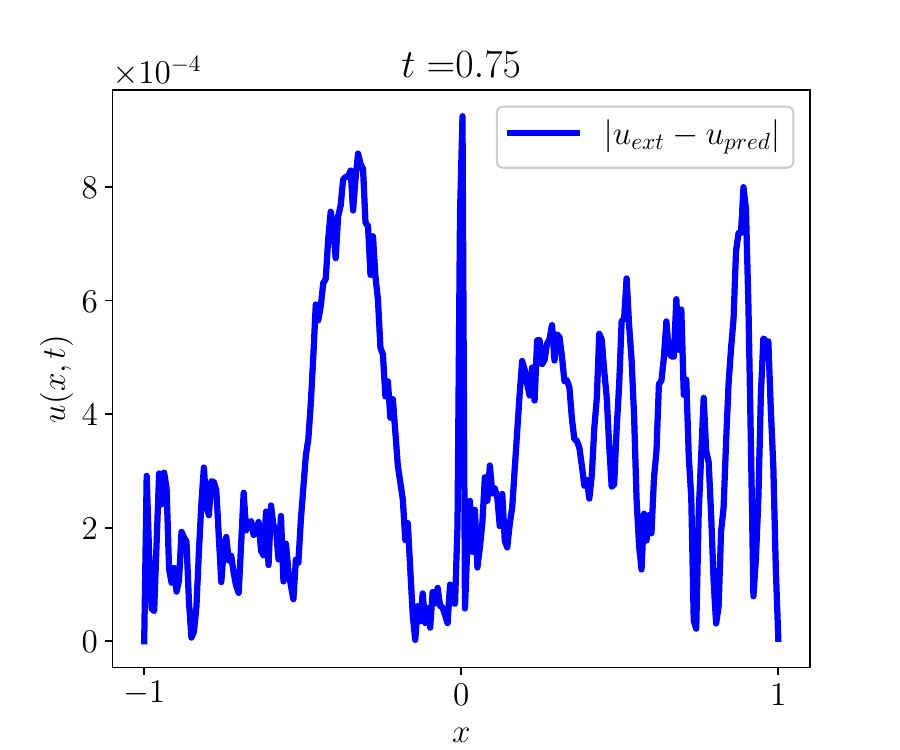}
	}	
	\subfigure[FFT of PIMWNN solution.]{
        \includegraphics[width=2.in]{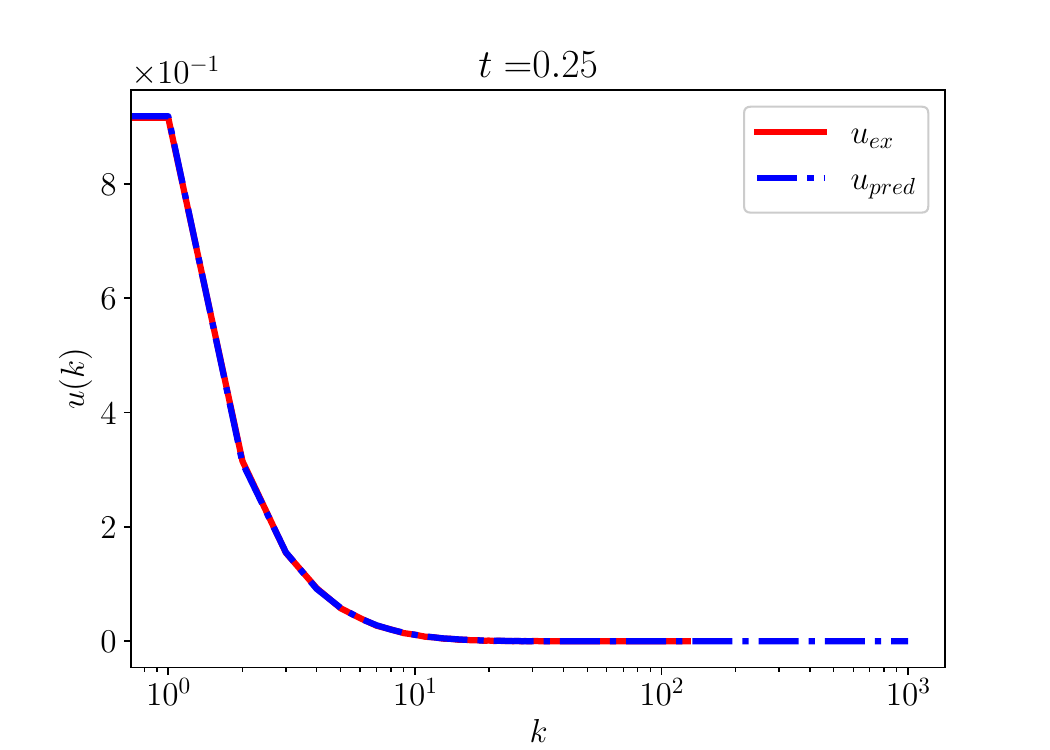}
        \includegraphics[width=2.in]{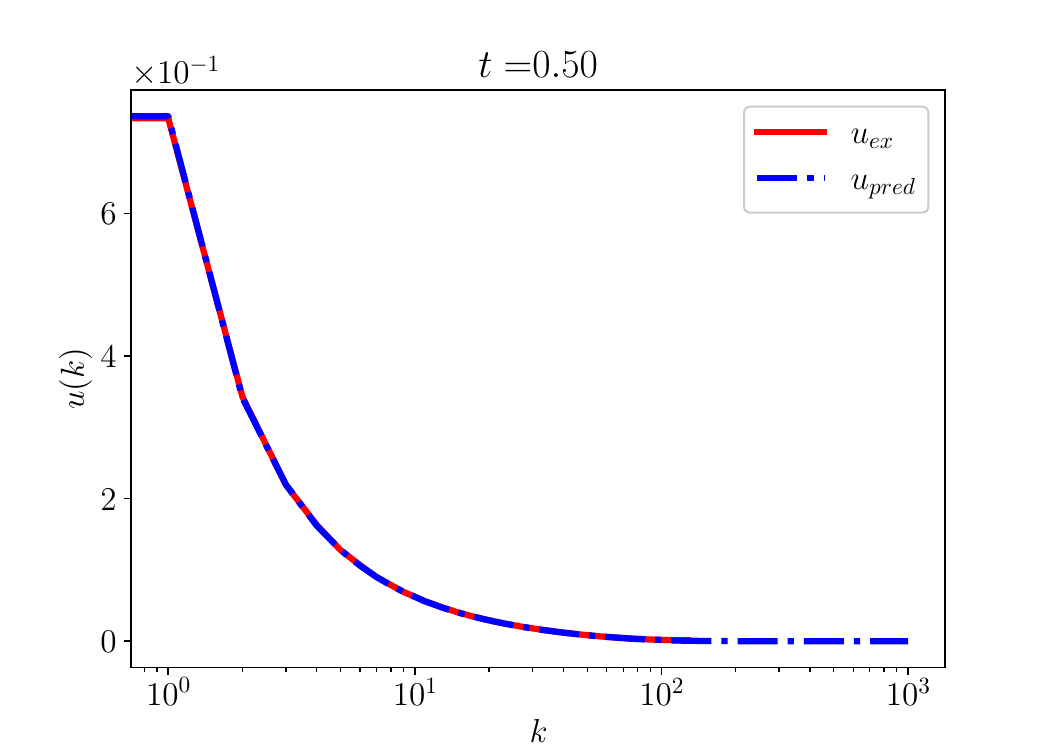}       
        \includegraphics[width=2.in]{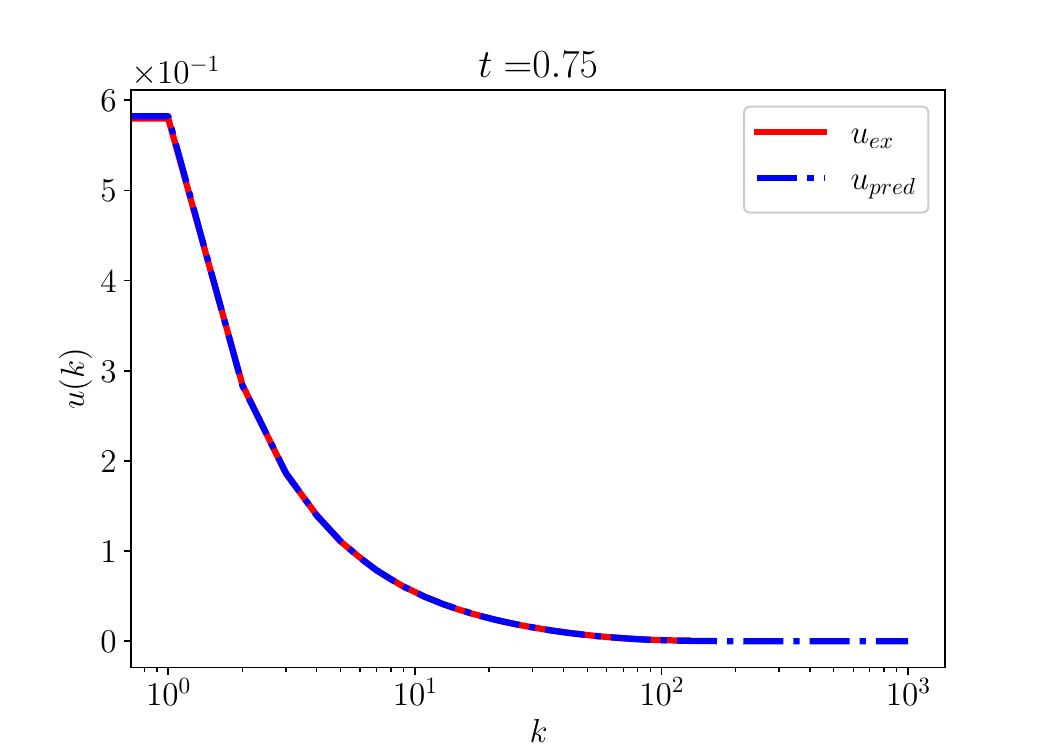}
	}	
\caption{Solution and error plots for Burger's equation at different time.}\label{tc-11}
\end{figure}

This non-linear equation is solved by the PIMWNN method as we mentioned in Section \ref{sc-bergers}. The parameters used in calculation are as follows: $dt = 0.001, N_f = 2000, J_0 = -1, J = 9$. In every time step, we use ten Picard iterations. Fig \ref{burgers-fig} shows the numerical results calculated by PIMWNN method and Randomized neural network method\cite{Wangfei2023}.  We can see that the solution accuracy of the proposed method is two orders of magnitude higher than that of the random neural network method. A comparison between the exact and the PIMWNN solutions at different time are shown in Fig. \ref{tc-11}(a). The pointwise error is shown in Fig. \ref{tc-11}(b). In every time step, the relative $L^2$ error are all of the order of $10^{-3}$. The FFTs of the PIMWNN solutions are shown in Fig. \ref{tc-11}(c). The results show that the spectrum still matches well with that of the exact solution after a long period of iteration.

\section{DISCUSSION}\label{discussion}

This paper aims to develop a physics-informed multiresolution wavelet neural network for solving partial differential equations. Based on multiresolution wavelet theory, we construct the multiresolution wavelet neural network. Combined with physics information which is proposed in PINN method, we propose the physics informed multiresolution wavelet neural network (PIMWNN) method. We apply the proposed method to solve both linear and nonlinear equations with diverse solution characteristics. The results demonstrate that the PIMWNN method surpasses PINNs in both efficiency and accuracy. The PIMWNN method offers several advantages: (i) The neural network employed in the PIMWNN method is meaningful, with its approximation capability strengthening as more scales are selected; (ii) It is mesh-free, and the algorithm's implementation is straightforward; (iii) Boundary conditions can be easily handled within this framework; (iv) The PIMWNN method efficiently solves linear time-dependent problems using a space-time approach; (v) It effectively captures high-frequency components of solutions and addresses the spectral bias issue in network training by increasing the number of wavelet scales.
On the other hand, in wavelet theory, high-dimensional wavelet bases are constructed by the direct product of 1D wavelet bases. For high-dimensional problems, this leads to a large number of wavelet bases, which translates into a large number of hidden nodes in the construction of wavelet neural networks. Meanwhile, this requires the least squares method to solve a linear system of equations with a large number of columns during network training. Although the above issues do not affect the correctness of the method in this paper, they do restrict its practicality and solution speed in actual applications. 

To address these issues, further research work can be carried out from the following aspects. First, conduct research on regional decomposition methods based on the approach in this paper. Decompose complex large-scale problems into several small and simple sub-problems for solution, thereby reducing the difficulty and complexity of problem-solving. In sub-regions, when the objective function is a relatively simple function, it can be approximated with a smaller scale, thus reducing the number of wavelet bases. This is also a research task planned for our next step. Second, combined with the characteristics of specific problems, carry out research on adaptive selection of wavelet bases (corresponding to adaptive adjustment of network structure). Third, conduct research on fast solution methods for least squares problems. We also plan to apply the PIMWNN method to a broader range of problems to explore its potential applications, including electromagnetic field problems, plasma physics problems, and fluid mechanics problems, among others.

\clearpage

\bibliographystyle{unsrt}

\bibliography{myreference}

\end{document}